\documentstyle{article}

\font\Bbb=msbm10 scaled \magstephalf

\def\codim{\mathop{\rm codim}}
\def\infp{\mathop{{\rm inf}\llap{\phantom{\rm p}}}}
\def\colon{:\;}
\def\<{\langle}
\def\>{\rangle}

\def\dst{\displaystyle}
\def\proof{\medbreak{\bf Proof.\ }}
\def\sbt{\subset}
\def\Pc{{\cal P}}
\def\R{\hbox{\Bbb\char'122}}
\def\x{{\bf x}}

\def\grad{\mathop{\rm grad}}
\def\pa{\partial}
\def\ob{{\bf 0}}

\def\vp{\varphi}
\def\la{\lambda}
\def\be{\beta}
\def\vola{\mathop{\rm vol}}
\def\al{\alpha}
\def\y{{\bf y}}

\def\ol{\overline}
\def\intt{\mathop{\rm int}}
\def\dsl{\displaylines}
\def\nomd{\noalign{\medskip}}
\def\vb{{\bf v}}
\def\ey{\emptyset}
\def\bcapi{\mathop{\hbox{$\bigcap$}}\limits}
\def\bcupi{\mathop{\hbox{$\bigcup$}}\limits}
\def\uli{\mathop{\hbox{$\underline{\lim}$}}\limits}
\def\oli{\mathop{\hbox{$\overline{\lim}$}}\limits}

\def\Kca{{\cal K}}
\def\ts{\times}
\def\de{\delta}
\def\iy{\infty}
\def\Pc{{\cal P}}
\def\diam{\mathop{\rm diam}}
\def\spt{\supset}

\def\ve{\varepsilon}
\def\N{\hbox{\Bbb\char'116}}
\def\eb{{\bf e}}
\def\Pc{{\cal P}}
\def\Kc{{\cal K}}
\def\La{\Lambda}
\def\ext{\mathop{\rm ext}}
\def\sm{\setminus}
\def\ev{\equiv}
\def\ab{{\bf a}}
\def\bb{{\bf b}}
\def\Lc{{\cal L}}
\def\ga{\gamma}
\def\si{\sigma}
\def\om{\omega}

\def\yb{{\bf y}}
\def\zb{{\bf z}}
\def\vr{\varrho}
\def\ub{{\bf u}}
\def\nosm{\noalign{\smallskip}}
\def\dist{\mathop{\rm dist}}
\def\CK{\overline{C}}

\def\dim{\mathop{\rm dim}}
\def\tr{\Delta}

\def\con{\mathop{\rm con}}

\def\cb{{\bf c}}

\def\grad{\mathop{\rm grad}}
\def\w{{\bf w}}

\begin{document}

\title
{THE GENERALIZED MINKOWSKI FUNCTIONAL
WITH APPLICATIONS IN APPROXIMATION THEORY}
\author
{SZIL\'ARD GY. R\'EV\'ESZ \thanks
{This joint research has been supported by a ``Marie
Curie Research Fellowship" of the European Community programme
``Improving Human Potential" under contract\break
\# HPMF--CT--2000--00670. The first author was supported in part
by the Hungarian National Foundation for Scientific Research,
Grant \# T034531 and T032872.}
\,\, and YANNIS SARANTOPOULOS}

\date{}

\maketitle

{\small\noindent{\bf Abstract.}}
We give a systematic and thorough study of geometric notions and results
connected to Minkowski's measure of symmetry and the extension of the
well-known Minkowski functional to arbitrary, not necessarily symmetric
convex bodies $K$ on any (real) normed space $X$.
Although many of the notions and results we treat in this paper
can be found elsewhere in the literature, they are scattered and
possibly hard to find. Further, we are not aware of a systematic
study of this kind and we feel that several features, connections
and properties
-- e.g. the connections between many equivalent formulations --
are new, more general and they are put in a better perspective now.
In particular, we  prove a number of fundamental properties of the
extended Minkowski functional $\alpha (K,\x)$, including
convexity,
global Lipschitz boundedness, linear growth and approximation of
the
classical Minkowski functional of the central symmetrization of
the body $K$.
Our aim is to present how in the recent years these notions
proved to be surprisingly relevant and effective in problems of
approximation theory.

\medskip
{\bf MSC 2000 Subject Classification:}{\it Primary 46B20;\\
Secondary 41A17, 41A63, 41A44, 26D05.}

{\bf Key words and phrases:\/} {\it convex body, support function,
supporting hyperplanes, halfspaces and layers,
Minkowski functional, convex functions in normed spaces,
Lipschitz bounds, central symmetrization, centroid, cone of
convex bodies, measure of symmetry, width of $K$ in a direction,
homothetic transformations, existence of minima of continuous
convex functions in normed spaces, separation of convex sets,
multivariate polynomials, Bernstein and Chebyshev type
extremal problems for multivariate polynomials.}

\vfil\eject

\section{Introduction}
The present work grew out of our studies
on inequalities and extremal problems for
multivariate polynomials, especially multivariate polynomials
defined on infinite dimensional normed spaces. For a general
introduction to polynomials on infinite dimensional spaces we
refer to the comprehensive first chapter of Dineen's
monograph [D]; as for the discussion of the particular
multivariate polynomial extremal problems connected to
generalized Minkowski functional, see \S9 of this article.

A quite general question in these type of problems presents itself
in the following setting. Suppose we are given a (continuous)
polynomial $p$ of
known degree $n$ and satisfying some boundedness condition on a
reasonably well-shaped set, say on a convex body $K\sbt X$, where
$X$ is a normed space. If $\Pc_n(X)$ is the space of polynomials
$p\colon X\to\R$ of degree not exceeding $n$, the ``initial
conditions" of a usual extremal problem are $p\in \Pc_n(X)$ and
${\|p\|}_K\le 1$, where ${\| p\|}_K:= \sup_K |p|$. The aim is to
find best possible (or at least of a correct order of magnitude)
estimates on some expressions --- say, a linear functional --- of
$p$, which depend on a point $\x\in X $. For example, the
``Bernstein problem" is to estimate $\|\grad p(\x)\|$; the ``point
value Chebyshev problem" is to estimate $\big|p(\x)\big|$. To
solve these questions, apart from functional analysis and
approximation theory, a crucial role is played by the geometric
features of the configuration of $K$ and $\x$. In particular, any
such estimate has to use some quantitative formulation of the
location of $\x$ with respect to the body $K$. Obvious choices
--- like the metric distance of $\x\notin K$ from $K$, or $\x\in
K$ from the boundary $\pa K$
--- yield only essentially trivial but too rough estimates.

If the {\it convex body\/} $K\sbt X $
(i.e., $K$ is a bounded, convex and closed set with nonempty interior)
is also {\it centrally symmetric\/} with respect to the origin, then
there is a norm ${\|\cdot\|}_{(K)}$ on $ X $ --- equivalent to the
original norm $\|\cdot\|$ of $ X $ in view of $B_{ X,\,
\|\cdot\|}(\ob,r)\sbt K\sbt B_{X,\,\|\cdot\|}(\ob,R)$ --- which
can be used successfully in approximation theoretic questions as
well. For this norm ${\|\cdot\|}_{(K)}$ the unit ball of $X$ will
be $K$ itself, $B_{X,\,{\|\cdot\|}_{(K)}}(\ob,1)= K $, while for
any $\x\in X$ the {\it Minkowski functional\/} or {\it (Minkowski)
distance function\/} [HC, p.~57] or {\it gauge\/} [Ro, p.~28] or
{\it Minkowski gauge functional\/} [P, 1.1(d)]
$$
\vp_K(\x):= \inf\{\la> 0\colon \x\in\la K\}
\leqno(1.1)
$$
will serve as definition of ${\|\x\|}_{(K)}$. Clearly (1.1) is a
norm on $ X $ if and only if the convex body $K$ is centrally
symmetric
with respect to the origin.

In case $K$ is nonsymmetric, (1.1) still can be used. But then,
even the choice of the homothetic centre is questionable since the
use of any alternative gauge functional
$$
\vp_{K, {\x_0}}(\x):= \inf\big\{\la> 0\colon \x\in
\la(K-\x_0)\big\}
$$
is equally well justified. Moreover, neither is good enough for
the applications.

There are two possible approaches here. The first one, which to us
seems to be more easily understood, is to define some
``continuous'', increasing, absorbing set function $t \rightarrow
K_t$, and to derive the corresponding generalized functional
analogously to (1.1). That will be pursued in this Introduction.
On the other hand, another general approach is described e.g. by
Gr\"unbaum in connection to the construction of measures of
symmetry in general, see [G, p. 245]. In this approach, one
defines first a functional $f(K,\x)$ by directly using some
geometric features of the configuration of $K$ and $\x$. Then a
measure of symmetry is defined by considering the extremal value
of $f$ over the whole of $K$ (or $X$). Depending on the geometric
properties we start with, corresponding approaches lead to
valuable information. In \S 4 we examine how this second approach
can be used to arrive at equivalent forms and thus unifying a
number of seemingly different geometric constructions.

An interesting construction of the first type, which resembles
somewhat to our upcoming choice of subject, is the so called
``floating body" construction \footnote{We are grateful to
Professor Apostolos Giannopoulos for calling our attention to
this notion and providing us the corresponding references.} (see
[BL]), which was successfully applied in a number of geometric
problems (c.f.\ e.g.\ the recent survey article [B\'a]).

Skipping the precise details of formulation, let us describe
heuristically this definition. If $0< \be < 1$, imagine that the
closed surface of $\pa K$ is a container which contains a liquid
of volume $\be \vola (K)/2$. Rotating the body shell (in the
fixed gravitational field) to all possible positions, the liquid
flows freely around. Then the body will consist of points which
became ``wet" during the rotation, and of ``dry'' points which
went under water in neither position. Cutting down the wet
points, the set of the dry points form a smaller set $K(\be)$
which is easily seen to be convex, and, naturally, form a
monotone set sequence of the continuous parameter $\be$.

The crucial feature is that, although $K$ is not symmetric, a
natural, geometrical definition is given to a reasonable set
sequence $K(\be)$, which reflects the geometric characteristics
(in our case the boundary shape in particular). If we want,
similarly to (1.1), a ``floating distance" can be defined:
$$
\psi_K(\x):= \inf\big\{\la > 0\colon \x\in K(1-\la)\big\},
$$
although this is not usual in the applications of the floating
body construction.

Now this is very much like what we want, but at least two
properties constitute some major drawback. First, the definition
is restricted to points $\x\in K $, and second, the reference to
the volume excludes infinite dimensional spaces. Moreover, we do
not know any particular applications in approximation theory
which could be well described using the floating body notion.

It is worth noting that for $\x\notin K$ there is an alternative
notion, the ``illumination body" construction, quite similar in
nature. See for instance [We]. However, this also relies on use of
the volume, hence is restricted to finite dimensional spaces.

So, our choice will be another construction of a set function
which was considered first by Hammer [H].
The use of an equivalent geometric notion
(the somewhat implicit use of $\la(K,\x)$, see (4.2) below)
has already appeared in approximation theory in the work of
Rivlin and Shapiro [RiSh] dealing with the
Chebyshev problem for finite dimensional spaces.
Although Rivlin--Shapiro, as well as Kro\'o--Schmidt [KS] later on,
have worked with exterior points of $K$, it was realized in [KR] that
the same definition works for points of $K$ as well. Also, in [KR] a few
observations and questions were mentioned concerning the generalized Minkowski
functional $\al(K,\x)$ (see definition (2.13) below).

However, to the best of our knowledge this present work is the
first attempt
in bringing together what geometers on the one hand, and
approximation
theorists on the other hand have contributed in the subject. For this reason,
new results
are accompanied by quite a number of facts which were already
known to either
geometers or approximation theorists, but not in our setting and
without the connections our investigation reveals. Lacking the drive
of accommodating to application needs of approximation theory, some
geometric facts, which were accessible to experts well before, seem to
surface or to be explored the first time here.
In particular, not even the fundamental properties of the
generalized
Minkowski functional, like convexity, Lipschitz boundedness, or
linear speed
growth, seem to have been established before. Moreover, in the
present work
we investigate the subject in the full generality of normed
spaces.

A careful checking reveals that essentially $\al(K,\x)$ does
not depend  on the initial norm of $X$, at least in the
sense that any equivalent norm can be considered. Actually the
sets $K^\la$ of (2.12), and hence $\al(K,\x)$, depend only on
geometry and not on the norm. Therefore, this notion can be equally
well introduced in topological vector spaces (t.v.s) where the
(topological) dual $X^*$, the support function and every necessary
ingredient is available. However, if $X$ is only a t.v.s, then the
mere existence of the convex body $K$ leads to a topologically
equivalent norm, the unit ball of which is the central
symmetrization $C$ (c.f. (2.1)) of $K$. As we will show in the
sequel, actually ${\|\cdot\|}_{(C)}= \vp_C$ and $\al (K,\cdot)$
have only a very small --- globally bounded --- difference, and so
the seemingly easy generalization to t.v.s. is essentially void.

A word of warning concerns the use of the inner product space
structure, which is justified in $\R^d$ in view of the equivalence
of all norms, but is not available for normed spaces in general.
Hence, here we must carefully distinguish the space $X$, its unit
ball $B=B_X(\ob,1)$ and its unit sphere $S=S_X(\ob,1)=\pa B$ on
one hand, and $X^*$ its corresponding unit ball $B^*$ and its unit
sphere $S^*$ on the other hand. The value of a linear functional
$\x^*\in X^*$ at a vector $\x\in X$ will be denoted by
$\<\x^*,\x\>$.

\section {Notations, definitions, and general background}

The {\it central symmetrization \/} of a set $K$ in a normed space
$X$ is
$$
C:= C(K):= {1\over 2}(K-K):= \left\{{1\over
2}(\x-\y)\colon \x,\y\in K \right\} .
\leqno(2.1)
$$
The central symmetrization of $K$ is centrally symmetric with
respect to the origin. In case $K$ is a convex body, we also have
$\ob\in \intt C$ \footnote{Throughout the paper we denote when
convenient
$C(K),\, \tau (K,\vb),\, \al (K,\x ), \,w(K,\vb^*)$, etc. by
$C , \, \tau ,\, \al , \, w$, etc., respectively.}.
On the other hand, even though $K$ is assumed
to be closed, $C$ is not necessarily closed (c.f.\ \S6 of [RS]),
hence $C$ is not a convex body in general. Nevertheless, the
closure $\ol C$ of $C$ is a symmetric convex body, which is also
{\it fat}, and $\intt C\sbt C\sbt \ol{\intt C}= \ol C$.

The {\it `` maximal chord} " of $K$ in direction of $\vb\ne \ob$
is
$$
\dsl{ \tau (K,\vb):= \sup\{ \la \ge 0:\exists \,\, \y,\,\zb\in K
\,\hbox{ s.t. }\zb=\y+\la\vb\}= \cr\nomd (2.2)\hfill =\sup \big\{
\la \ge 0\colon K\cap(K+\la\vb)\ne \ey\big\}= \hfill\qquad
\cr\nomd =\sup\{\la \ge 0\colon \la\vb\in K-K\}= 2\sup\big\{
\la>0\colon \la\vb \in C\big\}= \cr\nomd =2\max\big\{ \la \ge
0\colon \la\vb \in \ol{C}\big\}= \tau\big( C,\vb\big). \cr\nomd}
$$

Usually $\tau(K,\vb)$ is not a ``maximal" chord length, but only a
supremum, however we shall use the familiar finite dimensional
terminology
(see for example [W]).

The {\it support function\/} to $K$, where $K$ can be an arbitrary
set, is defined for all $\vb^*\in X^*$ (sometimes only for
$\vb^*\in S^*$) as
$$
h(K,\vb^*):= \sup_K\vb^*=\sup\big\{\<\vb^*,\x\>\colon \x\in
K\big\},
\leqno(2.3)
$$
and the {\it width of $K$ in direction\/} $\vb^*\in X^*$ (or
$\vb^*\in S^*$) is
$$
\dsl{(2.4)\hfill w(K,\vb^*):=
h(K,\vb^*)+h(K,-\vb^*)=\sup_K\vb^*+\sup_K(-\vb^*)= \hfill \qquad
\cr\nomd
=\sup\big\{\<\vb^*,\x-\y\>\colon \x,\y\in
K\big\}=2h\big(C,\vb^*\big)=w\big(C,\vb^*\big).
\cr}
$$
Let us introduce the notations
$$
X_t(\vb^*):= \big\{\x\in X \colon \<\vb^*,\x\>\le t\big\},
\qquad
X(K,\vb^*):=X_{h(K,\vb^*)}(\vb^*).
\leqno(2.5)
$$
Clearly the closed halfspace $X(K,\vb^*)$ contains $K$ and the
hyperplane
$$
H(K,\vb^*):=H_{h(K,\vb^*)}(\vb^*),
\quad
H_t(\vb^*):= \big\{\x\in X \colon \<\vb^*,\x\>= t\big\}=\pa
X_t(\vb^*)
\leqno(2.6)
$$
is a supporting hyperplane \footnote{Note that throughout the
paper we
mean ``supporting" in the weak sense, that is, we do not require
$K\cap H(K,\vb^*)\ne\ey$, but only $\dist
\big(K,H(K,\vb^*)\big)=0$.
The same convention is in effect for other supporting objects
as halfspaces, layers etc.} to $K$.

A {\it layer\/} (sometimes also called {\it strip\/} in the
literature) is the region of $X$
enclosed by two parallel
hyperplanes, i.e.
$$
L_{r,s}(\vb^*):= \big\{\x\in X \colon r\le \<\vb^*,\x\>\le s\big\}
=X_s(\vb^*)\cap X_{-r}(-\vb^*),
\leqno(2.7)
$$
while the {\it supporting layer\/} or {\it fitting layer\/}
of $K$ with normal $\vb^*$ is
$$
\dsl{(2.8)\hfill L(K,\vb^*):=
X(K,\vb^*)\cap X(K,-\vb^*)=L_{-h(K,-\vb^*),\,h(K,\vb^*)}(\vb^*)=
\hfill\qquad \cr\nomd
=\big\{\x\in X \colon -h(K,-\vb^*)\le
\<\vb^*,\x\>\le h(K,\vb^*)\big\}.
\cr}
$$
By convexity, $K$ is the intersection of its `` supporting
halfpaces " $X(K,\vb^*)$, and grouping opposite normals we get
$$
K=\bcapi_{\vb^*\in S^*}X(K,\vb^*)=\bcapi_{\vb^*\in S^*}
L(K,\vb^*).
\leqno(2.9)
$$
Any layer (2.7) can be homothetically dilated with quotient $\la
\ge 0$ at any of its symmetry centers lying on the symmetry
hyperplane
$H_{r+s\over2}(\vb^*)$
to obtain
$$
L^\la_{r,s}(\vb^*):= \left\{\x\in X \colon {\la+1\over
2}r-{\la-1\over 2}s\le \<\vb^*,\x\>\le {\la+1\over
2}s-{\la-1\over 2}r\right\}.
\leqno(2.10)
$$
In particular, we have also defined
$$
L^\la(K,\vb^*)=L_{-h(K,-\vb^*),\,h(K,\vb^*)}^\la (\vb^*)
\leqno(2.11)
$$
and by using (2.11) one can even define
$$
K^\la:= \bcapi_{\vb^*\in S^*} L^\la(K,\vb^*) .
\leqno(2.12)
$$

Note that $\;K^\la$ can be empty for small values of $\la$. Using
the convex, closed, bounded, increasing and (as easily seen, c.f.
Proposition 3.3) even absorbing set system ${\{K^\la\}}_{\la\ge
0}$, the {\it generalized Minkowski functional\/} or {\it gauge
functional\/} is defined as
$$
\al(K,\x):= \inf\{\la\ge 0\colon \x\in K^\la\}.
\leqno(2.13)
$$

Denoting the cone of convex bodies in $ X $ by $ \Kc :=
 \Kc ( X )$, in  analogy to e.g.\ [HC, Chapter 1.9, p.~302], we see
that $\al$ is a mapping from $\Kca \ts X \to \R_+$. Speaking of
topological or analytical features of $\al$, we always understand
$\Kca$ as the topological space equipped with the
(Pompeiu--)Hausdorff(--Blaschke) metric $\de=\de^H$ (c.f.\ [HC]
p.~307)
$$
\de(K,M):= \max\left\{
\sup_{\x\in K} \infp_{\y\in M} \|\x-\y\|,
\ \sup_{\y\in M} \infp_{\x\in K}\|\x-\y\|\right\},
\leqno(2.14)
$$
while $\R_+=[0,\iy)$ is considered with its usual metric. Note
that (2.14) may be defined for any subsets of $X$. We shall also use
the notation
$$
\al_K:= \inf \al (K,\cdot)=\inf\{\la\ge 0\colon K^\la\ne \ey\}.
\leqno(2.15)
$$

For further use, we introduce the notation
$$
\dsl{(2.16)\hfill \Pc\colon \Kca\ts \R_+\to \Kca\hfill\qquad
\cr\nomd
\Pc(K,\la):= \Pc^\la(K):= K^\la
\cr}
$$
for the formation of the `` $\la$-powered set " of $K$. To avoid
notational complications, here and elsewhere if necessary, we
identify $\Kca$ and its (metric) complete closure of not
necessarily proper (i.e., the set of all convex, closed, bounded,
nonempty, but not necessarily $\intt K\ne\ey$) convex bodies.
Cauchy
sequences in $\Kca$ converge.

In reference to the notations (2.1)--(2.4), we also use the
notions of the {\it width of\/} $K$, the {\it diameter of\/} $K$
and the {\it inradius\/} of $C$ which are defined as
$$
\dsl{(2.17)\hfill w:= w(K):= \inf\big\{w(K,\vb^*)\colon
\vb^*\in S^*\big\},\hfill\qquad
\cr\nomd
(2.18)\hfill d:= d(K):= \diam\, (K):=
\sup\big\{\|\x-\y\|\colon \x,\y\in K\big\},\hfill\qquad
\cr}
$$
and
$$
r:= r\big(C):= \sup\big\{r\colon B(\ob,r)\sbt C\big\},
\leqno(2.19)
$$
respectively.

Finally, we introduce the notation
$$
D:= D(K):= \sup_{\x\in K}\|\x\|.
\leqno(2.20)
$$

In view of (2.1)--(2.4) we have
$$
w(K)=w\big(C\big),
\leqno(2.21)
$$
and by taking into consideration the notation (2.20), we also have
$$
d(K)=d\big(C\big)=2 D\big(C\big).
\leqno(2.22)
$$
Observe also that $C \sbt \ol{C}$ implies
$r(C) \le r(\ol{C})$. Moreover, we claim that
$$
r\big(C\big)=\max \big\{r\colon B(\ob,r)\sbt \ol{C}\big\}=
r(\ol{C}).
\leqno(2.23)
$$
Indeed, the maximum exists because the ball $
B\big(\ob,r\big(\ol{C}\big)\big)$
is the closure of the union of all
balls in $\ol{C}$ and obviously it is contained in
$\ol{C}$ by closedness. But $r(\ol{C})$
can not exceed $r\big(C\big)$
since by fatness $\intt B\big(\ob,r\big(\ol{C}\big)\big)\sbt \intt
\ol{C}\sbt C$. Whence
$r\big(\ol{C}\big)=r\big(C\big)$ proving (2.23).

Usually maximal chord length and width are considered in finite
dimensional spaces, equipped with the inner product structure,
and so they can be compared (after identifying $X$ and $X^*$).
 See e.g.\ [W], Th. 7.6.1. For infinite dimensional spaces this is
not possible, but some of the intrinsic connections can be proved.

\proclaim {Lemma 2.1}.
$w(K)=2r\big(C\big)$.

\proof
In view of (2.21) is enough to show that $w(C)=2r(C)$ . Also
by fatness we may assume that $C$ is closed, that is $\ol C = C$.
So, for $C= \ol C$ take now $B_r=B\big(\ob,r(C)\big)\sbt C$
(c.f.\ (2.23)). In view of (2.4) and (2.17) we have to prove that
$r(C)=\inf\big\{h(C,\vb^*)\colon \vb^*\in S^*\big\}$. The ``$\le$ "
direction
is obvious as $B_r\sbt C$, $\|\vb^*\|=\sup_{B} \vb^*=
{1\over r}\sup_{B_r} \vb^*$ and (2.3) yields
$r(C)=\sup_{B_r} \vb^*\le h(C,\vb^*)$, for all $\vb^* \in
S^*$.

On the other hand, by using the definition (2.19) of $r(C)$,
for any $\ve >0$ there exists some point $\x\notin C$ with
$\|\x\|\le r+\ve$. Applying the separation theorem to this
$\x$ and $C$, we conclude that $r+ \ve\ge \|\x\|\ge \<\vb^* ,\x\>\ge
h(C,\vb^*)$, for some $\vb^*  \in S^*$. Taking infimum we conclude
also the `` $\ge$ " part of the assertion.

\proclaim {Corollary 2.2}.
$w(K)=\inf\big\{\tau(K,\vb)\colon \vb\in S\big\}$.

\proof
Since by (2.2) $\pa C=\left\{{\tau(K,\vb)\over 2}\vb \colon \vb\in
S\right\}$, the inradius $r\big(C\big)$ is the infimum of
$\tau(K,\vb)/2$ over all $\vb \in S$. Therefore a
reference to Lemma 2.1 proves the statement.

\medskip
Note that again by the definition of $\|\vb^*\|$, for any
$\vb^*\in S^*$, one can always choose a $\vb\in S$ so that
$\<\vb^*,\vb\>> 1-\ve$. Hence, the obvious relation
$$
w(K,\vb^*) \ge \tau(K,\vb) |\<\vb^*,\vb\>| \leqno(2.24)
$$
also implies the `` $\ge$ " part of Corollary 2.2, while
the other direction can be also obtained  from the following
result.

\proclaim {Lemma 2.3}.
For any $\vb\in S$ there exists $\vb^*\in
S^*$ such that
$$
w(K,\vb^*)=\tau (K,\vb)\cdot \<\vb^*,\vb\>.
\leqno(2.25)
$$

\proof
Consider a normal vector $\vb^*\in S^*$ to
$C$ at ${\tau\over 2}\vb\in \pa C$.
Clearly $h(C,\vb^*)= {\tau\over 2}\<\vb^*,\vb\>$.
Now (2.25) follows from the fact that $w(C,\vb^*)= w(K,\vb^*)$
(see (2.4)).

\medskip
At this stage we can formulate the generalization of [W] Theorem
7.6.1 to infinite dimensional spaces.

\proclaim {Proposition 2.4}. For all $\vb \in S$, $\vb^*\in S^*$
we have the following sharp estimates
$$
w(K)\le \tau(K,\vb)\le d(K),\quad w(K)\le w(K,\vb^*)\le d(K)\,.
$$

\proof The lower bounds follow from Corollary 2.2 and the
definition (2.17), respectively. The upper bound $\tau(K,\vb)\le
d(K)$ is an easy consequence of definitions (2.18) and (2.2).
These three inequalities are obviously sharp.

Now consider inequality (2.24) and take the supremum over $S^*$ on
both sides. Since $\sup\big\{\<\vb^*,\vb\>\colon \vb^*\in
S^*\big\}=\|\vb\|= 1$, we get
$\tau(K,\vb)\le\sup\big\{w(K,\vb^*)\colon \vb^*\in S^*\big\}$. In
view of $d(K)=\sup\big\{\tau(K,\vb): \,\vb \in S\big\}$, we have
$d(K)\le\sup\big\{w(K,\vb^*):\, \vb^*\in S^*\big\}$. To prove the
converse inequality, by (2.2)--(2.4) is enough to consider $C$
in place of $K$. But this is trivial since
$w(C,\vb^*)=2h(C,\vb^*)= 2\sup_{\pa C} \vb^*= 2\sup
\left\{\left\<\vb^*,{\tau\over 2}\vb \right\>:\, \vb\in
S\right\}\le 2\sup_S \tau/2= d(K)$.

\proclaim Example 2.5. \rm Let
$\ell^2$ be the real Hilbert space of square summable sequences
and let
$$
 K:= \left\{\x={(x_n)}_{n=1}^\iy \in \ell^2\colon
\sum\limits_{n=1}^\iy w_n x_n^2\le 1\right\}\,,
$$
where the weights $w_n$ satisfy $1< w_n<2$ $(n\in\N)$. Then $K$ is
a closed, bounded, convex set in $\ell^2$ with $\ob$ in its
interior. However,
\begin{itemize}
\item[(i)]
with $w_n:= 1+{1\over n}$ $(n\in\N),$ we have $d(K)=2$, but the
supremum in (2.18) is not attained;
\item[(ii)]
with $w_n:= 2-{1\over n}$ $(n\in \N)$, we have $w(K)=\sqrt 2$,
but the infimum in (2.17) is not attained.
\end{itemize}

\proof By definition, $K$ is closed. Since $B(\ob, 1/\sqrt 2)
\sbt K\sbt B(\ob,1)$, we have $\sqrt 2\le w(K)\le d(K)\le 2$. Also
$\ob\in \intt B(\ob,1/\sqrt 2)\sbt \intt K$ and $K$ is obviously
convex. Moreover, $K$ is a centrally symmetric convex body and
therefore $w(K)=2\inf \big\{\|\x\|\colon \x\in\pa K\big\}$ and
$d(K)=2\sup \big\{\|\x\|\colon \x\in \pa K\big\}$. Let
${(\eb_n)}_{n=1}^\iy$ be the standard unit vector basis of
$\ell^2$ where $\eb_n:= (0,\dots,0,1,0,\dots,0,\dots)$, with $1$
at the $n^{\rm th}$ place. In case (i), ${1\over \sqrt{w_n}}
\eb_n\in\pa K$, hence $d(K)\ge {2\over \sqrt{w_n}}$ $(n\in \N)$
and we find $d(K)=2$. On the other hand, for $w_n>1$ $(n\in\N)$ it
is clear that any $\x\in K\setminus\{\ob \}$ satisfies
$\sum\limits_n
x_n^2 <\sum\limits_n w_n x_n^2\le 1$, hence $\|\x\|<1$, and $d(K)$
is not attained.

Similarly, in case (ii), for any $\vb^*=\vb=(v_n)\in
S_{\ell_2}^*=S_{\ell_2}$ we get
$$
h(K,\vb^*)\ge \left\< \vb^*, {1\over \sqrt{\sum\limits_n w_n
v_n^2}} \vb \right\> ={1\over \sqrt
{\sum\limits_n w_n v_n^2}}>{1\over \sqrt 2
{\|\vb^*\|}_2}={1\over \sqrt 2},
$$
and hence $w(K, \vb^*)>\sqrt 2$. Now since $w(K, \eb_n)=2\left\|
{1\over \sqrt{w_n}}\cdot \eb_n\right\|={2\over \sqrt{w_n}}\to \sqrt 2$,
we have $w(K) =\sqrt 2$. But $w(K)=\sqrt 2$ is not attained.

\section{Basic properties}
\proclaim {Lemma 3.1}.
The mapping $\Pc$, defined in $(2.16)$, has the
following monotonicity properties:
\begin{itemize}
\item[i)]
For $K\in \Kc$ fixed, $\Pc ^\la(K)=K^\la$ is increasing
on $\R_+$.
\item[ii)]
For $\la\ge 1$ fixed, $\Pc ^\la(K)=K^\la$ is increasing
on $\Kc$.
\item[iii)]
$\Pc$ is increasing on $\Kc\ts[1,\iy)$.
\item[iv)]
For $\la<1$ fixed, $\Pc ^\la(K)$ is not monotonic.
\end{itemize}

\proof
i) By definition (2.12), it suffices to show the
assertion for an arbitrary layer $L$ which is obvious.\\
ii) Again, it suffices to show the assertion for layers. Let
$L\sbt L'$; then they have to be parallel and in view of (2.7)we
must have $L=L_{r,s}(\vb^*)$, $L'=L_{r',s'}(\vb^*)$ and
$[r,s]\sbt [r',s']$, for some $\vb^*\in S^*$ and $r,s,r',s'\in
\R$. Now (2.10) obviously gives $L^\la\sbt {L'}^\la$, as stated.
Note that $\vb^*$ brings the question for layers to a question
for real intervals. Moreover, it is easy to see e.g.\ from (2.10)
that actually for any layer $L$ and any $\la\ge 0$ we have
$$
L^\la={\la+1\over 2}L+{\la-1\over 2}(-L)=\la\cdot C(L)+\x_0
\hbox{ with } L-\x_0=-(L-\x_0)
\leqno(3.1)
$$
i.e. with any symmetry center $\x_0$ of $L$. (In fact (3.1) is
valid for any centrally symmetric convex set).\\
(iii) Follows from i) and ii).\\
iv) Let $L\sbt L'$ with common normal vector $\vb^*$,
say. If $r+s=r'+s'$, (where $L=L_{r,s}(\vb^*)$,
$L'=L'_{r',s'}(\vb^*)$), then $L^\la\sbt {L'}^\la$ for any
$\la\ge 0$. But if  $\la<1$, $r=r'=0$, $s=1$ and $s'>{1+\la\over
1-\la}$,
then $L^\la$ and ${L'}^\la$ are disjoint.

\proclaim {Lemma 3.2}.
$K^\la=\bcapi_{\mu>\la}K^\mu$.

\proof We have
$\bcapi_{\mu>\la}K^\mu=\bcapi_{\mu>\la}\bcapi_{\vb^*\in
S^*} L^\mu(K,\vb^*)=\bcapi_{\vb^*\in S^*} \bcapi _{\mu>\la}
L^\mu(K,\vb^*)=\bcapi _{\vb^*\in S^*} L^\la(K,\vb^*)=K^\la$, as
stated.

\proclaim {Proposition 3.3}. $\{ K^\la\}_{\la\ge 0}$ is absorbing,
and thus $\al(K,\x)$ is finitely defined all over $ X $. Moreover,
$\al(K,\x)$ is bounded on bounded sets.

\proof Since $\intt K\ne \ey$, there exists a ball $B(\x_0,r)\sbt
K$. By ii) of Lemma 3.1, for any $\la\ge 1$ we have $K^\la\spt
B^\la(\x_0,r)=B(\x_0,\la r)$, and thus $\bcupi_{\la\ge 0}
K^\la\spt \bcupi_{\la\ge 1}B(\x_0, \la r)= X $.  Boundedness of
$\al(K,\x)$ on bounded sets follows from
$$
B(\ob,R)\sbt B(\x_0,\la r)\sbt K^\la \quad \bigg(\la\ge
\max\Big\{1,{1\over r}\big(\|\x_0\|+R\big)\Big\}\bigg)\,.
$$

\proclaim {Corollary 3.4}.
For all $\x\in X $ we have $\x\in
K^{\al(K,\x)}$, and in fact (2.13) is a minimum:
$$
\al(K,\x)=\min\{\la\ge 0\colon \x\in K^\la\}.
\leqno(3.2)
$$

\proof
By Lemma 3.1 i), Proposition 3.3 implies that $\La:=
\{\la\ge 0\colon \x\in K^\la\}$ is a nonempty interval with right
endpoint at $\iy$. By Lemma 3.2, $\La$ is closed, too.

\proclaim {Proposition 3.5}.
$\al(K,\x)\to\iy$ as $\|\x\|\to\iy$.

\proof Since $K$ is bounded, $K\sbt B(\ob,R)$ for some $0<R<\iy$,
whence for any $\la\ge 1$ by Lemma 3.1 ii) we have $K^\la\sbt
B(\ob, \la R)$. It follows that for any $N\in \N$ and $\|\x\|>NR$
the point $\x\notin K^N$ and thus $\al(K,\x)\ge N$ $\big(\|\x\| >
R_0(N):= NR)$.

\proclaim {Proposition 3.6}.
For fixed $K\in\Kc$, $\al(K,\cdot)$ is
a convex function on the whole space $ X $.

\proof
Let us denote $\al':=\al(K,\x')$, $\al'':=
\al(K,\x'')$, $\be:= t\al'+(1-t)\al''$ and
$\hat\x:=t\x'+(1-t)\x''$, $0\le t\le 1$.
Moreover, for a given
$\vb^*\in S^*$ denote $L:=L(K, \vb^*)=L_{r,s}(\vb^*)$, with $r:=
-h(K,-\vb^*)$, $s:= h(K,\vb^*)$.

It is easy to see, e.g.\ by using  (3.1), that
$$
L^\be(K,\vb^*)=L^{t\al'+(1-t)\al''}=tL^{\al'}+(1-t)
L^{\al''}.
\leqno(3.3)
$$
In view of Corollary 3.4 we have $\x'\in K^{\al'}\sbt L^{\al'}$,
$\x'' \in K^{\al''}\sbt L^{\al''}$ and therefore $\hat {\x}\in
tL^{\al'}+(1-t)L^{\al''}=L^\be(K,\vb^*)$. Since this is true for all
$\vb^*$, it follows that $\hat {\x} \in
\bcapi_{\vb^*\in S^*} L^\be(K,\vb^*)=K^\be$. Hence
$\al(K,\hat\x)\le \be$ and the Proposition is proved.

\proclaim {Corollary 3.7}.
For fixed $K\in \Kc$, $\al(K,\cdot)$ is
continuous on $ X $.

\proof If $ X $ is finite dimensional, say $ X =\R^d$,
convexity itself implies continuity, see for instance [Ro, Corollary
10.1.1.]. However, for infinite dimensional spaces finiteness and
convexity do not imply continuity. Thus, here we have to invoke
also local boundedness, furnished by the second half of
Proposition 3.3. Now, a reference to [P, Proposition 1.6 and the
following Remark] yields continuity; moreover, it follows that
$\al(K,\cdot)$ is {\it locally Lipschitzian.\/} Instead of
utilizing this stronger assertion, we postpone the question of
Lipschitz bounds to later analysis. We will prove even more
precise results in Proposition 5.3 and Theorem 5.5.

\proclaim {Corollary 3.8}.
{\rm (i)} $K^\la=\big\{\x\in X \colon \al(K,\x)\le\la\big\}$.
\begin{itemize}
\item[(ii)]
$\intt K^\la=\bcupi_{\mu<\la}K^\mu=\bcupi_{\mu<\la}\intt K^\mu$.
\item[(iii)]
$\al(K,\x)<\la \quad \hbox {if and only if } \quad \x\in \intt
K^\la$.
\item[(iv)]
$\al(K,\x)=\la\quad \hbox{if and only if }\quad \x\in \pa K^\la$.
\item[(v)]
$\al(K,\x)>\la \quad \hbox{if and only if }\quad \x\in\ext
K^\la:= X \sm K^\la$.
\item[(vi)]
$\al_K<1 \quad\hbox{for every}\quad K\in \Kc( X )$.
\end{itemize}

\proof
(i) The proof is obvious in view of (2.13) and Lemma 3.1 i).\\
(ii) If $\la>\mu$, by Lemma 3.1 i) $K^\la\spt K^\mu\spt \intt
K^\mu$. Hence $\intt K^\la\spt \intt K^\mu$
and the `` $\spt$ " inclusion follows. On the other hand,
for any $\x\in \intt K^\la$ there exists a ball $B(\x,r)\sbt
K^\la=\big\{\x\in X \colon \al(K,\x)\le \la\big\}$ in view of i).
If the
value of $\al(K,\x)$ were $\la$, then $\al$ would have a local
maximum at $\x$ which, in view of convexity, would imply that
$\al$ is a constant multiple of $\la$. But this contradicts
Proposition 3.5. Thus
$\al(K,\x)<\la$, i.e.\ $\x\in \bcupi_{\mu<\la}K^\mu$, and the
first equality follows. Equality of the middle and the right
expression is obvious from the first identity.\\
(iii) follows from (i) and (ii), and (iv) follows from (i)
and (iii), while (v) follows from (i) trivially. Finally, (ii) and
and the fact that $\intt K^1=\intt K\ne \ey$ imply (vi).

\section{Relatives of $\al(K,\x)$}
There are various other geometric quantities, which are closely
related to our definition of the Minkowski functional. Many of
these equivalent forms and related quantities play a crucial
role in approximation theoretic applications. Therefore, we find
it appropriate to give a more or less comprehensive account of
them here. The use of them will also be rather helpful in the
forthcoming sections.

Let us start with the description of the original definition of
Minkowski,
[Mi], as presented in Gr\"unbaum's article, c.f. [G p. 246].
Denote
$$
t:= t(K,\vb^*,\x):=
{2\<\vb^*,\x\>-h(K,\vb^*)+h(K,-\vb^*)\over w(K,\vb^*)}.
\leqno(4.1)
$$

For fixed $\vb^*$ this function is an affine linear functional in
$\x\in X$,
while for fixed $\x$ it is a norm-continuous
mapping from $S^*$ (or $X^*\setminus \{\ob \}$) to $\R$.
In fact, for fixed
$\vb^*\in S^*$, $t$ maps the layer $L(K,\vb^*)$ to
$[-1,1]$, and $L^\eta(K,\vb^*)$ to $[-\eta,\eta]$. Therefore, the
two forms of the following definition are really equivalent;
$$
\dsl{(4.2)\hfill \la:= \la(K,\x):=
\sup \big\{\eta>0\colon \exists \vb^*\in S^*,\;\x\in \pa
L^\eta(K,\vb^*)\big\}
\hfill\qquad
\cr\nomd
=\sup\Big\{\big|t(K,\vb^*,\x)\big|\colon \vb^*\in S^*\Big\}
=\sup\Big\{ t(K,\vb^*,\x) \colon \vb^*\in S^*\Big\}
.\cr}
$$
Note that $t(K,\vb^*,\x)=-t(K,-\vb^*,\x)$ and therefore we don't
have to use the absolute value.

In fact, $\la(K,\x)$ expresses the supremum of the ratios
of the distances between the
point $\x$ and the symmetry hyperplane $\frac12
(H(K,\vb^*)+H(K,-\vb^*))$
of any layer $L(K,\vb^*)$ and the
half-width $w(K,\vb^*)/2$. Now Minkowski's definition was
$$
\varphi(K,\x):= \inf \left\{
\frac{\min\{\dist\left(\x,H(K,\vb^*)\right), \dist\left(\x, H(K,-\vb^*)\right)\}}
{\max\{\dist\left(\x, H(K,\vb^*)\right), \dist\left(\x,H(K,-\vb^*)\right)\}}
\, : \, \vb^*\in S^* \right\},
$$
which clearly implies the relation
$$
\varphi(K,\x)=\frac{1-\la(K,\x)}{1+\la(K,\x)} \qquad (\x \in K).
$$
Although this $\varphi(K,\x)$ seems to be used traditionally
only for $\x\in K$, extending
the definition to arbitrary $\x \in X$ yields the similar relation
$$
\varphi(K,\x)=\frac{|1-\la(K,\x)|}{1+\la(K,\x)} \qquad (\x \in X).
$$

\proclaim {Proposition 4.1}.
We have $\al(K,\x)=\la(K,\x)$
$(\forall \x\in X)$.

\proof
If $\x\in \pa L^\la(K,\vb^*)$ for some $\la>0$ and
$\vb^*\in S^*$, then obviously
$\x\notin L^{\la-\ve}(K,\vb^*)$
for any $\ve>0$, and hence also $\x\notin K^{\la-\ve}$ and
$\al(K,\x)\ge \la$. This proves $\al(K,\x)\ge \la(K,\x)$. Let
now $\vb^*\in S^*$ be arbitrary. Then, in view of (4.2) we have
$\x\in \pa L^{|t|}(K,\vb^*)\sbt L^\la(K,\vb^*)$ with $t$ and
$\la$ defined in (4.1) and (4.2), respectively. Therefore $\x\in
K^\la$ and thus $\al\le \la$ follows.
\medskip

We now define some
geometric quantities using homothetic images. Let us define
$$
\dsl{(4.3)\hfill \be(K,\x):= \sup\big\{\la\ge
0\colon \{-\la(K-\x)+\x\} \sbt K\big\}\qquad(\x\in K), \hfill\quad
\cr\nomd
(4.4) \hfill \vr(K,\x):= \sup \Big\{0<\la<1\colon K\cap
\big\{\x+\la(K-\x)\big\}=\ey\Big\}\quad (\x\notin K),\hfill
\cr\nomd
(4.5) \hfill \vr^*(K,\x):=
\sup\Big\{0<\la<1 \!\colon\! \intt K\cap
\big\{\x+\la(K-\x)\big\}=\ey\Big\} \quad \hfill (\x\notin K).
\cr}
$$

Let us mention that these geometric quantities were used in various
ways for
long. In particular, Schneider defines the set of parallel bodies
of a
body $K$ with respect to another one, say $D$, c.f. [Sch p.~134] or
(7.15) below. When $D=K$,
this reduces to parallel bodies of $K$ closely related to the
above $\beta$.
In fact, in [Note 3.1.14, p. 141 Sch] Schneider mentions that
$r(K,-K)=\max_K\beta$ is Minkowski's measure of symmetry.
Moreover, the above
defined $\rho(K,\x)$ is closely connected to a
well-known notion, the so-called ``associated
bodies'' of Hammer, c.f. [Sch, pp. 141-142] or [H].
Thus it will be used to establish the connection
of the parametric set sequence $K^\la$ to Hammer's associated
bodies.

In the following we also consider the set of straight lines
$$
\dsl{(4.6) \hfill \Lc:= \Lc(K,\x):= \big\{l=\{\x+t\vb\colon t\in
\R\big\}:\#\{K\cap l\}>1, \ \vb\in S\big\}=\hfill\qquad
\cr\nomd
=\big\{l=\{\x+t\vb \colon t\in \R\}\colon \vb\in S,\; K\cap
l=[\ab,\bb], \ \ab\ne \bb\big\}.
\cr}
$$
For $\x\in\intt K$, $\Lc$ contains all straight lines through $\x$,
but for $\x\notin K$ only some lines are included. When they need
to be distinguished, we always assume that the notation of
the endpoints of $K\cap l$ is chosen so that $\|\x-\bb\|\le
\|\x-\ab\|$. With these in mind, we define
$$
\dsl{(4.7)\hfill \mu(K,\x):=
\inf\left\{{\|2\x-\ab-\bb\|\over \|\bb-\ab\|}\colon l\in
\Lc\right\}
\quad (\x\notin K ),\hfill\qquad
\cr\nomd
\nu(K,\x):= \sup\left\{{\|2\x-\ab-\bb\|\over
\|\bb-\ab\|}\colon l\in \Lc\right\}\quad (\x\in K ),
\cr\nomd
\si(K,\x):= \inf\left\{{\|\bb-\x\|\over \|\ab-\x\|}\colon \; l\in
\Lc\right\} \quad (\x\in X ),
\cr\nomd
\om(K,\x):= \sup\left\{{\|\ab-\x\|\over \|\ab-\bb\|}\colon \; l\in
\Lc\right\}\quad (\x\in K ),
\cr\nomd
\ga^2(K,\x):= \inf\left\{{\|\bb-\x\|\,\|\ab-\x\|\over
\left\|{1\over 2}(\bb-\ab)\right\|^2}\colon \;l\in
\Lc\right\}\quad
(\x\in K).
\cr}
$$
Note the various restrictions on the domains of the definition
of these quantities. Note also that in fact $\si(K,\x)$ is exactly
the
equivalent form of $\varphi(K,\x)$ mentioned also in [G p.~247].
This form
is used by even more authors than the original one, hence we could
refer to
common knowledge for $\varphi(K,\x)=\si(K,\x)$. Still, we give
proofs for
all equivalences below.

\proclaim {Proposition 4.2}.
i) In $(4.5)$, the sup is actually a maximum.
\begin{itemize}
\item[ii)]
$\vr(K,\x)=\vr^*(K,\x)$.
\item[iii)]
There exists $\vb^*\in S^*$ such that the supporting
hyperplane (2.6) separates (weakly) $K$ and $\x+\vr(K-\x)$,
in other words
$$
h(K,\vb^*)\le (1-\vr)\<\vb^*,\x\>-\vr h(K,-\vb^*).
\leqno(4.8)
$$
\end{itemize}

\proof
i) Consider for some scalar $0<\si <1$ the set
$
M:= M_\si:= \intt K\cap\big(\x+\si (K-\x)\big).
$
Let us suppose first that $M \ne \ey$. As $\x+\si(K-\x)$ is
also a fat set, this entails that $\intt M\ne \ey$, and so
there exist $\yb\in M$ and $r>0$ such that even $B(\yb,r)\sbt
M$. On the other hand, boundedness of $K$ implies $K\sbt B(\x,
R)$, for $R$ large enough. Take now $\ve<r\si/R$ and
consider $\zb:= {\si-\ve\over \si}(\yb-\x)+\x$. Clearly
$\|\zb-\yb\|={\ve\over \si}\|\yb-\x\|<r$, hence even $\zb\in
B(\yb,r)\sbt \intt K$. Since
$\yb\in\x+\si(K-\x)$, we have
$\zb\in\x+(\si-\ve)(K-\x)$ and this shows that together with
$M_\si\ne\ey$ also $M_{\si-\ve}\ne\ey$ follows with some
$\ve>0$. Thus $M_{\vr^*}\ne\ey$ would lead to a contradiction
with (4.5), and so $M_{\vr^*}=\ey$, as needed.\\
ii) According to part i) the convex sets $\intt K$ and
$\x+\vr^*(K-\x)$ are disjoint. Applying the separation theorem of
convex sets (c.f.\ e.g.\ [Ru], Theorem 3.4.) we conclude the
existence of a $\vb^*\in S^*$ separating them. Since $\ol{\intt
K}=K$, $\vb^*$ separates even $K$ and $\x+\vr^*(K-\x)$ at least
in the weak (i.e. $\le$) sense. However, this entails that $K$
and $\x+(\vr^*-\ve)(K-\x)$ are separated by $\vb^*$ even in the
strong sense for any $\ve>0$. Hence, in particular $K\cap
\big\{(\vr^*-\ve)(K-\x)+\x\big\}=\ey$, and $\vr\ge \vr^*-\ve$
$(\forall \ve>0)$, which leads to $\vr\ge \vr^*$. On the other
hand, it is clear that for all $\x\notin K$ we have $0<\vr \le
\vr^*<1$ and thus we arrive at $\vr=\vr^*$.\\
iii) Since $\rho = \rho^*$, by construction we see that
$\vb^*$ separates $K$ and $\x+\vr(K-\x)$, and a simple calculation
gives (4.8)
as well.

\proclaim {Proposition 4.3}. For all $\x\notin K$ we have
$\si(K,\x)=\vr(K,\x)$.

\proof Let $\la<\si:=\si(K,\x)$ be arbitrary. Then by definition
(4.7) of $\sigma, \,\, K$ and $\x+\la(K-\x)$ are disjoint, hence
also $\la \le \vr$, which implies $\si\le \vr$.

On the other hand, in view of Proposition 4.2 iii) there exists
$\vb^*$ separating $K$ and $\x+\vr(K-\x)$, i.e. $K-\x$ and
$\vr(K-\x)$, which implies even for arbitrary $\ab,\bb\in K$
(and not only for points on the same line $l$) that with some
constant $c>0$
$$
\<\vb^*,\bb-\x\>\ge c\ge \big\<\vb^*,\vr(\ab-\x)\big\>\ge 0.
\leqno(4.9)
$$
Clearly, separation occurs only if $\<\vb^*, \yb-\x\>$ does not
vanish for points
$\yb\in\intt K$ at least. However, this yields
$\|\bb-\x\|\ge \vr\|\ab-\x\|$ from (4.9), and by definition
(4.7) we get even $\si\ge \vr$ after taking the infimum over
$l\in\Lc$. Together with the first part this completes the proof.

\proclaim {Proposition 4.4}.
For all $\x\notin K$ we have
$\si(K,\x)={\mu(K,\x)-1\over \mu(K,\x)+1}$.

\proof
For any $\x\notin K$ it is clear that $\mu>1$. Moreover,
$$
\dsl{{\mu-1\over \mu+1}=1-{2\over \mu+1}=1-{2\over
\inf\limits_\Lc
\big(\|2\x-2\ab\|/\|\bb-\ab\|\big)}=
\cr\nomd
=1-\sup\limits_\Lc{\|\bb-\ab\|\over
\|\x-\ab\|}
=\inf\limits_\Lc{\|\bb-\x\|\over \|\ab-\x\|}=\si,
\cr}
$$
as asserted.

\proclaim {Proposition 4.5}.
For all $\x\in X\sm K$ we have
$\vr(K,\x)={\al(K,\x)-1\over \al(K,\x)+1}$.

\proof
Suppose first that $\x\in K^\la=\bcapi_{\ub^*\in S^*}
L^\la(K,\ub^*)$. Then a simple calculation leads to the
inequality
$$
\<\ub^*,\x\>\le {\la\over 2}w(K,\ub^*)+{1\over
2}\big(h(K,\ub^*)-h(K,-\ub^*)\big)\quad (\forall \ub^*\in S^*)
$$
which is inherited by $\al=\inf\{\la\colon \x\in K^\la\}$ as well.
Thus, also for the particular $\vb^*$, furnished by Proposition
4.2 iii), we obtain
$$
\<\vb^*,\x\>\le {\al\over 2}w(K,\vb^*)+{1\over
2}\big(h(K,\vb^*)-h(K,-\vb^*)\big).
\leqno(4.10)
$$
On combining (4.8) and (4.10) we are led to
$$
h(K,\vb^*)\le (1-\vr)\left\{{\al+1\over 2}
h(K,\vb^*)+{\al-1\over 2} h(K,-\vb^*)\right\}-\vr h(K,-\vb^*)
$$
and some calculation yields
$$
w(K,\vb^*)\big((\al-1)-\vr(\al+1)\big)\ge 0.
\leqno(4.11)
$$
Since $w(K,\vb^*)\ge w (K)>0$, (4.11) entails $\vr\le {\al-1\over
\al+1}$.

Let now $1<\la<\al$, i.e.\ $\x\notin L^\la(K,\ub^*)$ for some
appropriate $\ub^*\in S^*$. Then, a simple calculation shows that
already $\x+{\la-1\over \la+1}\big(L(K,\ub^*)-\x)\big)$ is
disjoint from $L(K,\ub^*)$, thus even more so
$\left\{\x+{\la-1\over \la+1}(K-\x)\right\}\cap K=\ey$, and
${\la-1\over \la+1}\le \vr$ follows.  Letting here $\la\to \al^{-}$
we even get ${\al-1\over \al+1}\le\vr$, and the opposite
inequality being already proven, this concludes the proof.

\proclaim {Proposition 4.6}.
For $\x\in K$ we have $\si(K,\x)= {1\over \omega(K,\x)}-1$.

\proof
$$
\si=\inf_\Lc {||\bb-\x||\over ||\ab-\x||}= \inf_\Lc
{||\bb-\ab||\over
||\ab-\x||}-1 = {1\over \sup\limits_\Lc {||\ab-\x||\over
||\bb-\ab||}}
-1= {1\over \omega} -1.
$$

\proclaim {Proposition 4.7}.
For $\x\in K$ we have $\si(K,\x)=\be(K,\x)$.

\proof
Let $\ab\in\pa K$ be any boundary point. It is easy to
see (by closedness of $K$) that (4.3) is actually a maximum,
that is, $-\be(K-\x)+\x\sbt K$, and therefore $\ab':=
-\be(\ab-\x)+\x\in K$. Consider now an $l\in \Lc$, $K\cap
l=[\ab,\bb]$, $\|\bb-\x\|\le \|\ab-\x\|$, and estimate
${\|\bb-\x\|\over \|\ab-\x\|}\ge {\|\ab'-\x\|\over
\|\ab-\x\|}=\be$. We only need to take infimum over $l \in \Lc$
to obtain $\si \ge \be$.

For the other direction let us consider the set
$K':= -\si(K-\x)+\x$. Let $\yb\in K$ be arbitrary: we show
that also $\zb:= -\si(\yb-\x)+\x\in K$, leading to $K'\sbt K$
and $\si\le \be$, which will prove Proposition 4.7.

Take now the straight line $l$ through $\x$ and $\yb$ $(\yb\ne
\x$, as for $\yb=\x$ also $\zb=\x$ and there remains nothing to
prove). Let $\ab$ be the endpoint of $K\cap l$ on the side of
$\yb$, --- i.e.\ $\yb\in[\x,\ab]$ and $\ab\in\pa_{rel} K\cap l$.
Then by
definition (4.7) we must have for the other endpoint $\bb\in K\cap
l$, on the other side of $\x$, the inequality $\|\bb-\x\|\ge
\si\|\ab-\x\|\ge \si\|\yb-\x\|=\|\zb-\x\|$. This shows that
$\zb\in [\x,\bb]\sbt K$ and hence $\zb\in K$ as needed.

\proclaim {Proposition 4.8}.
For $\x\in K$ we have
$\be(K,\x)={1-\la(K,\x)\over 1+\la(K,\x)}$.

\proof For arbitrary $\vb^*\in S^*$ consider the layer $L\colon =
L(K, \vb^*)$. By (4.1) a simple calculation gives $\{-{1-|t|\over
1+|t|}(L-\x)+\x\} \sbt L$ $\big(t:= t(K,\vb^*,\x)\big)$. Since the
function $\vp(t):= {1-t\over 1+t}$ $(0\le t\le 1)$ is strictly
decreasing, this last relation is also true for
$\la=\sup_{S^*}|t|$ defined by (4.2). Thus $-{1-\la\over
1+\la}(L-\x)+\x\sbt L$ $(\vb^*\in S^*)$, and hence
$$
\dsl{(4.12) \hfill K=\bcapi_{\vb^*\in S^*} L(K,\vb^*)\spt
\bcapi_{\vb^*\in
S^*}\left\{-{1-\la\over
1+\la}\big(L(K,\vb^*)-\x\big)+\x\right\}\hfill\qquad
\cr\nomd
=-{1-\la\over
1+\la}(K-\x)+\x,
\cr}
$$
proving that $\be\ge {1-\la\over 1+\la}$.

To deduce the converse, consider now $-\be(K-\x)+\x\sbt K$.
Forming supporting layers normal to $\vb^*\in S^*$, we find
$$
-\be\big(L(K,\vb^*)-\x\big)+\x=
L\Big(\big(-\be(K-\x)+\x\big),\vb^*\Big)\sbt
L(K,\vb^*),
$$
that is, $-\be(L-\x)+\x\sbt L$ $(\forall \vb^*\in S^*)$. For any
given $\vb^*\in S^*$, however, the maximal scalar to satisfy
this relation is exactly ${1-|t|\over 1+|t|}$, hence $\be\le
{1-|t|\over 1+|t|}$ $(\vb^*\in S^*)$, and
taking infimum on the right hand side now gives $\be\le
{1-\la\over 1+\la}$, referring to the decreasing shape of
$\vp(t)$ again.

\proclaim {Proposition 4.9}.
For arbitrary $\x\in K$ we have
$\nu(K,\x)={1-\si(K,\x)\over 1+\si(K,\x)}$ .

\proof
Taking care of the distinction between $\ab$ and $\bb$ in
the notation for $K\cap l=[\ab,\bb]$ prescribed for $l\in
\Lc$, we write with $\vp(t):={1-t\over 1+t}$
$$
\dsl{\nu=\sup_{\Lc}{\|2\x-\ab-\bb\|\over
\|\bb-\ab\|}=\sup_{\Lc}{\|\x-\ab\|-\|\bb-\x\|\over
\|\x-\ab\|+\|\bb-\x\|}=
\cr\nomd
=\sup_{\Lc}{1-{\|\bb-\x\|\over \|\ab-\x\|}\over
1+{\|\bb-\x\|\over \|\ab-\x\|}}=
\sup_{\Lc} \varphi \Big( {\|\bb-\x\|\over \|\ab-\x\|} \Big)
={1-\inf_{\Lc}{\|\bb-\x\|\over
\|\ab-\x\|}\over 1+\inf_{\Lc}{\|\bb-\x\||\over
\|\ab-\x\|}}={1-\si\over 1+\si}\,.
\cr}
$$

\proclaim {Proposition 4.10}.
For $\x\in K$ we have
$\ga^2(K,\x)=4\omega(K,\x)\left(1-\omega(K,\x)\right)$.

\proof
Working in a similar fashion as in the proof of Proposition 4.9,
but using here
also the decreasing shape of $\psi(s):= 4s(1-s)$ on
$[{1\over2},1]$,
we can write
$$
\dsl{\ga^2=\inf_{\Lc}{\|\x-\ab\|\cdot\|\bb-\x\|\over
\left({1\over
2}\|\bb-\ab\|\right)^2}=\inf_{\Lc}\left\{4(1-s)s\colon
{1 \over 2}\le s={\|\ab-\x\|\over\|\bb-\ab\|}\le 1\right\}=
\cr\nomd
=\psi\left(\sup_{\Lc}{\|\ab-\x\|\over
\|\ab-\bb\|}\right)=\psi\left(\omega\right)=
4\omega(1-\omega).
\cr}
$$

\proclaim {Corollary 4.11}. We have the following relations for
$\al(K,\x)$ $:$
$$
\al(K,\x)=\la(K,\x)=\left\{
\begin{array}{ll}
\dst{1-\si(K,\x)\over 1+\si(K,\x)}=\dst{1-\varphi(K,\x)\over
1+\varphi(K,\x)}
& \quad (\x\in K)\\[4mm]
\dst{1+\si(K,\x)\over 1-\si(K,\x)}=\dst{1+\varphi(K,\x)\over
1-\varphi(K,\x)}
& \quad (\x\notin K)
\end{array} \right.
\hbox{for}\quad \x\in X,
\leqno{\rm (i)}
$$
$$
\dsl{{\rm(ii)} \hfill\al(K,\x)=\nu(K,\x)={1-\be(K,\x)\over
1+\be(K,\x)}=\sqrt{1-\ga^2(K,\x)}
= 2\om(K,\x) -1\hfill\qquad
\cr\nomd
\qquad \qquad \qquad \qquad \qquad \qquad \qquad \qquad \qquad
\qquad \qquad
\hbox{\qquad for} \quad \x\in K,
\cr\nomd
{\rm (iii)} \hfill\al(K,\x)=\mu(K,\x)={1+\vr(K,\x)\over
1-\vr(K,\x)}={1+\vr^*(K,\x)\over 1-\vr^*(K,\x)}
\qquad \hbox{\quad for}\quad \x\notin K.\hfill
\cr}
$$

\section{Lipschitz bounds and linear growth}
\proclaim {Lemma 5.1}. Let $L$ be any layer and let $0<\be<\ga$ be
arbitrary. Then we have
$$
\dist\{\pa L^\be,\pa L^\ga\}={\ga-\be\over 2}w(L).
$$

\proof This distance formula follows easily from (2.4) and (3.1).

\proclaim {Lemma 5.2}. Let $\x\ne \yb\in X $ and let $\vb:=
{1\over \|\x-\yb\|}(\x-\yb)$. Then
$$
\big|\al(K,\x)-\al(K,\yb)\big|\le {2\|\x-\yb\|\over
\tau(K,\vb)}\,, \leqno{\rm (5.1)}
$$
with $ \tau(K,\vb)$ defined by $(2.2)$. Moreover, $(5.1)$ is best
possible.

\proof For an arbitrary $\vb^*\in S^*$ let $L:= L(K,\vb^*)$. There
exist unique scalars $\be, \ga\ge 0$, so that $\x\in \pa L^\be$
and $\yb\in\pa L^\ga$. If we take $\de:= \|\x-\yb\|$, then it is
obvious that $\dist\{\pa L^\be,\pa L^\ga\}=\de\cdot
|\<\vb^*,\vb\>|$. On the other hand, by Lemma 5.1 this equals to
${|\ga-\be|\over 2}w(L)={|\ga-\be|\over 2}w(K,\vb^*)$. Hence, in
view of (2.24) we find
$$
\ga \le \be+{2\de|\<\vb^*,\vb\>|\over w(K,\vb^*)}
\le \be+{2\de \over \tau(K,\vb)}\,.
\leqno{\rm (5.2)}
$$
Now if $\x\in \pa L^\be$
and $\al(K,\x)=\al$, then by definition (4.2) and Proposition 4.1
we must have $\be\le \al$. But then (5.2) implies $\ga\le
\al+2\de/\tau$, whence $\yb\in L^{\al+2\de/\tau}$ for all
$\vb^*\in S^*$ and so $\y \in K^{\al+2\de/\tau}$ and
$\al(K,\yb)\le \al+2\de/\tau$. By symmetry we also have $\al\le
\al(K,\yb)+2\de/\tau$. These last two inequalities yield (5.1).

To show that (5.1) is sharp, consider $\x=(s+\de)\vb$, $\yb=s\vb$
(with $s$ large enough). For the given vector $\vb\in S$, take now a
$\vb^*\in S^*$ furnished by Lemma 2.3 satisfying (2.25). Consider
the point $\y =s\vb$, the affine linear functional (4.1) and the
quantity (4.2). By Proposition 4.1 we can write
$$
f(s):=\al(K,s\vb)=\la(K,s\vb)\ge t(K,\vb^*, s\vb). \leqno(5.3)
$$
Now, by taking into consideration (4.1) and (2.25), (5.3)  easily leads to
$$
\lim\limits_{s\to\iy}f(s)/s\ge {2\<\vb^*,\vb\>\over
w(K,\vb^*)}={2\over
\tau(K,\vb)}.
\leqno(5.4)
$$
Since the function $f(s)= \al(K,s\vb)$ is convex on $\R$, we get
$$
\dsl{(5.5)\hfill \sup\{\al(K,\x)-\al(K,\y): \x, \y \in X,
\x-\y=\de\vb \} \hfill \qquad
\cr\nomd
\ge \limsup\limits_{s\to
\iy}\big(f(s+\de)-f(s)\big) \ge \de\lim\limits_{s\to \iy} f(s)/s
\,.
\cr}
$$
It follows from (5.4) and (5.5) that (5.1) is sharp.

\proclaim {Proposition 5.3}.
We have the uniform Lipschitz bound
$$
\big|\al(K,\x)-\al(K,\yb)\big|\le {2\over w(K)}\|\x-\yb\|\quad
(\x,\yb \in X ).
\leqno(5.6)
$$

\proof
Combining Lemma 5.2 and Proposition 2.4 yields (5.6).

\proclaim {Lemma 5.4}.
Let $K,K'\in \Kc$, and $\x\in X$ be
arbitrary with the Hausdorff distance $\de=\de(K,K')<{1\over
2}w(K)$. Then we have
$$
\big|\al(K,\x) -\al(K',\x)\big|\le {\de(K,K')\big(4\|\x\|+4 D(K)+2
w(K)\big)\over w(K)\big(w(K)-2\de(K,K')\big)}.
\leqno(5.7)
$$

\proof
For an arbitrary $\vb^*\in S^*$ consider
the affine linear functional $t$ defined by (4.1) and similarly
the affine linear functional $t'$, where
$$
t':= t(K',\vb^*,\x):=
{2\<\vb^*,\x\>-h(K',\vb^* )+h(K',-\vb^*)\over w(K',\vb^*)}.
\leqno(5.8)
$$
By definition (2.14) of the Hausdorff distance
we have
$\big|h(K,\vb^*)-h(K',\vb^*)\big|
\break \le \de$,
$|w-w'|\le 2\de\,(<w)$, and thus
$$
\begin{array}{rl}
|t-t'|&\le \dst{2\big|\<\vb^*,\x\>\big| {2\de\over
ww'}+{2\de\over w}+ 2D(K') {2\de\over ww'}\le}\\[4mm]
&\le \dst{{4\de\big(\|\x\|+D(K)+\de\big)\over
w(w-2\de)}+{2\de(w-2\de)\over w(w-2\de)}=}\\[4mm]
&= \dst{\de\big(4\|\x\|+4D(K)+2w\big)\over w(w-2\de)}\,.
\end{array}
\leqno(5.9)
$$
In view of the definition (4.2), (5.9) implies
$$
\big|\la(K,\x)-\la(K',\x)\big|\le
{\de\big(4\|\x\|+4D(K)+2w\big)\over w(w-2\de)}\,.
\leqno(5.10)
$$
Now a reflection to Proposition 4.1 completes the proof
of the Lemma.

\proclaim {Theorem 5.5}.
Let $(K,\x), (M,\yb)\in\Kc\ts X $ with
$\de(K,M)<{1\over 2}w(K)$.
Then we have
$$
\dsl{(5.11)\hfill  \big|\al(K,\x)-\al(M,\yb)\big|\le
{\de(K,M)\big(4\|\yb\|+4D(K)+2w(K)\big)\over
w(K)\big(w(K)-2\de(K,M)\big)}
\hfill \qquad
\cr\nomd
+{2\|\x-\yb\|\over w(K)}.
\cr}
$$

\proof
On applying the triangle inequality
$\big|\al(K,\x)-\al(M,\yb)\big|\le\break
\big|\al(K,\x)-\al(K,\yb)\big|+\big|\al(K,\yb)-\al(M,\yb)\big|$,
we can
use Proposition 5.3 for the first term and Lemma 5.4 for the
second to obtain (5.11).

\medskip
Already in Proposition 3.5 we have seen
that $\al(K,\x)\to\iy$ $\big(\|\x\|\to\iy\big)$, and (5.1) with
(5.5)
suggest a linear speed growth. We can make this observation
much more precise.

\proclaim {Lemma 5.6}.
Denoting the classical Minkowski
functional (c.f.\ $(1.1)$, $(2.1)$) of $\ol{C}$ as
$\vp_{\ol{C}}={\|\cdot\|}_{\ol{C}}$, we have
$$
\vp_{\ol{C}}(\x)={2\|\x\|\over \tau\left(K, {\x\over
\|\x\|}\right)}\quad (\ob\ne \x\in X ).
\leqno(5.12)
$$

\proof
In view of (2.2) and the fatness of $C$ we see that
$\tau:=
\tau(K,\vb)=\tau\big(C,\vb\big)=\tau\big(\ol{C},\vb\big)$.
Choosing here $\vb\colon=\x/\|\x\|$, central symmetry of
$\ol{C}$
shows $\ol{C}\cap \{\la\vb\colon \la\in \R\}=\left[-{\tau\over
2}\vb,{\tau\over 2}\vb\right]$. Hence definition (1.1), applied
to $\ol{C}$ in place of $K$, gives
$\vp_{\ol{C}}(\x)={2\over \tau} \|\x\|$, as stated.

\proclaim {Theorem 5.7}.
For every $\x\in X$, $\x\neq\ob$, we have
$$
\left|\al(K,\x)-\vp_{\ol{C}}(\x)\right|=\left|\al(K,\x)-{2\|\x\|\over
\tau( K, \x/\|\x\|)}\right|\le {2D(K)\over w(K)}-1.
\leqno(5.13)
$$

\proof
Lemma 5.6 implies the equality part of (5.13). Now let
$\vb^*\in S^*$ be arbitrary and consider the affine
linear functional (4.1) for both $K$ and $\ol {C}$. As in
Lemma 5.6, we can easily see that $h\big(\ol{C},\vb^*\big)=
h\big(\ol{C},-\vb^*\big)$ by symmetry, and that
$w(K,\vb^*)=w\big(\ol{C},
\vb^*\big)$ by (2.4) and fatness of $C$. Therefore we find
$$
\dsl{(5.14)\hfill
t(K,\vb^*,\x)-t\big(\CK,\vb^*,\x\big)=
{-h(K,\vb^*)+h(K,-\vb^*)\over
w(K,\vb^*)}=\hfill\qquad
\cr\nomd
={2h(K,-\vb^*)-w(K,\vb^*)\over w(K,\vb^*)}=2{h(K,-\vb^*)\over
w(K,\vb^*)}-1\le {2D(K)\over w(K,\vb^*)}-1.
\cr}
$$
Taking into account the similar calculation
$$
t\big(\CK,\vb^*,\x\big)-t(K,\vb^*,\x)={2h(K,\vb^*)\over
w(K,\vb^*)}-1\le {2D(K)\over w(K,\vb^*)}-1,
$$
now we get by (2.17)
$$
\Big|t(K,\vb^*,\x)-t\big(\CK,\vb^*,\x\big)\Big|\le {2D(K)
\over w(K,\vb^*)}-1 \le {2D(K)\over w(K)}-1.
\leqno(5.15)
$$
By definition (4.2) we obtain from (5.15)
$$
\left|\la(K,\x)-\la\big(\CK,\x\big)\right|\le {2D(K)\over w(K)}-1.
\leqno(5.16)
$$
By using Proposition 4.1 and the fact that $\al\big(\CK,\x\big)=
\vp_{\CK}(\x)$, (5.16) gives (5.13).

\proclaim {Corollary 5.8}.
For $\x\to\iy$, $\al(K,\x)$ grows
linearly. In particular,
$$
\lim\limits_{s\to\iy}{\al(K,s\vb)\over s}={2\over
\tau(K,\vb)}\quad
(\vb\in X, \,\|\vb\|=1)
\leqno(5.17)
$$
uniformly on $S$, and also
$$
{2\over d(K)}=\uli_{\x\to\iy}{\al(K,\x)\over
\|\x\|}\le \oli_{\x\to\iy} {\al(K,\x)\over \|\x\|}={2\over
w(K)}.
\leqno(5.18)
$$

\proof
By the second term of (5.13), we can replace $\al(K,\x)$
by $2\|\x\|/\tau$,
$\big(\tau=\tau (K,\x/\|\x\|)\big)$ in (5.17)--(5.18). This gives
(5.17)
immediately, while (5.18) follows in view of the first part of
Proposition 2.4.

\proclaim {Corollary 5.9}. For the Hausdorff distance $\de$,
defined by (2.14), we have
\begin{itemize}
\item[i)]
$\de\big(K^\la,\la C)$ is bounded uniformly for $\la
> \al_K$ and
\item[ii)]
${1\over \la}K^\la \to C$ $(\la\to +\iy)$.
\end{itemize}

\proof
By convexity, $\de\big(K^\la,\la C\big)=\de\big(\pa K^\la,
\pa\la C\big)$. Applying Corollary 3.8 iv) and Lemma 5.6 we
obtain
$$
\dsl{(5.19)\hfill\de\big(K^\la,\la C\big)
=\de\big(\pa K^\la,\pa\la C\big)
\le\qquad\hfill\cr
\nomd
\hfill \le\sup\left\{
\left\|s\vb-\la{\tau(K,\vb)\over 2}\vb\right\|
\colon \vb\in S,\, 0<s=s(\la,\vb)
\quad\hbox{s.t.}\quad
\al(K,s\vb)=\la \right\}.\hfill\cr}
$$
For any particular $\vb$, $s(\la,\vb)$ exists provided the ray
$\{s\vb\colon s > 0\}$ intersects $K^\la$; this is certainly the
case when $\la>\al(K,\ob)$, since then $\ob\in\intt K^\la$, see
Corollary 3.8 iii). Hence, there exists some $s=s(\la,\vb)$ as
described in (5.19). Moreover, again by local boundedness (see
Proposition 3.3), $s\to\iy$ as $\la\to\iy$.

Here we refer to Theorem 5.7, using (5.13) with $\x=s\vb$ to
obtain
$$
\left|\al(K,s\vb)-{2s\over \tau(K,\vb)}\right|\le {2D(K)\over
w(K)}-1.
\leqno(5.20)
$$

Comparing (5.19) and (5.20), with $\la=\al(K,s\vb)$, we get
$$
\de\big(K^\la,\la C\big)\le \sup_{\vb\in S}{\tau(K,\vb)\over
2}\left({2D(K)\over w(K)}-1\right) < d(K)\cdot {D(K)\over w(K)},
\leqno(5.21)
$$
in view of Proposition 2.4. Clearly (5.21) implies i) and ii)
follows easily from i).

\section{Global properties on $\Kc\ts X $}
Perhaps the most important consequence of Theorem 5.5 is the
next result.

\proclaim {Corollary 6.1}. The generalized Minkowski functional
$\al$ is continuous on $\Kc\ts X $.

\proclaim {Proposition 6.2}. Let $K,M\in\Kc$. Then we have
$$
\al_{K+M}\le\max\{\al_K,\al_M\}\,,
\leqno(6.1)
$$
where $\al_{K+M}$, $\al_K$ and $\al_M$ are defined by $(2.15)$.

\proof
It is easy to see that for any $\vb^*\in S^*$
$$
L(K+M,\vb^*)=L(K,\vb^*)+L(M,\vb^*),
\leqno(6.2)
$$
hence also for any $\la\ge 0$
$$
L^\la(K+M,\vb^*)=L^\la(K,\vb^*)+L^\la(M,\vb^*)\,.
\leqno(6.3)
$$
But (6.3) implies that
$$
\dsl{
(6.4) \hfill
{(K+M)}^\la=\bcapi_{\vb^*\in S^*} L^\la(K+M,\vb^*)\spt\hfill\qquad
\cr\nomd
\spt \bcapi_{\vb^*\in S^*} L^\la (K,\vb^*)+\bcapi_{\vb^*\in S^*}
L^\la(M,\vb^*) = K^\la+M^\la.
\cr}
$$
Therefore, if $\yb\in K^\la$, $\zb\in M^\la$ and $\x=\yb+\zb$,
then $\x\in {(K+M)}^\la$ and thus $\al(K+M,\x)\le \la$ which
implies $\al_{K+M}\le \la$. As the argument applies to all
$\la>\max\{\al_K,\al_M\}$ (for then $K^\la\ne\ey$ and $M^\la\ne
\ey$), we can conclude (6.1).

\medskip
In the following we investigate the
case when $\al_K=0$. We will show that this characterizes
centrally symmetric domains.

\proclaim {Lemma 6.3}.
Suppose $\al_K=0$. Then
$K^{\al_K}=K^0= \{\x_0\}$, that is, $K^0$ is a one-point set.

\proof
We have
$$
\dsl{d(K^\la)=\sup_{\vb^*\in S^*}w(K^\la,\vb^*)\le \sup_{\vb^*\in
S^*} w(\big(L^\la(K,\vb^*),\vb^*\big)= \cr\nosm (6.5)\hfill
\hfill\qquad \cr\nosm =\sup_{\vb^*\in S^*} \la
w\big(L(K,\vb^*),\vb^*\big)=\la \sup_{\vb^*\in S^*}w(K,\vb^*)=\la
d(K), \cr}
$$
hence in case $\la\to 0^{+}$ we also have $d(K^\la)\to 0$. Since
$K^\la$ are closed sets and $K^\la$ $(\la\ge 0)$ forms a monotone
system, all conditions of a Cantor-type intersection theorem are
satisfied and by using also Lemma 3.2 we get
$$
K^0=\bcapi_{\la>0} K^\la\ne \ey.
$$

Since $d(K^\la)\to 0$, $|K^0|\le 1$ is obvious.

\proclaim {Question 6.4}. \rm For any $0\le \al_K<1$, is it true
in general that the minimal $K^{\al_K}$ of $K^\la$ is nonempty? In
other words, does $\al(K,\x)$ attain its infimum $\al_K$? Compare
also to Remarks 8.2.

In the next two theorems we recover the fact that Minkowski's
measure of symmetry is indeed a measure of symmetry in the sense
defined by Gr\"unbaum, see [G] or [HC p.~362]. That goes back to
Minkowski [Mi] and Radon [R] and is a well-known fact for finite
dimensional spaces. Note that already Gr\"unbaum points out, c.f.
[G, p.~234], that together with $f$ also $h\circ f$ must be
measure of symmetry for any homeomorphism $h$ of [0,1] onto itself
with $h(1)=1$. Thus, the equivalent formulations listed in
Corollary 4.11 lead to a measure of symmetry simultaneously and
the next two theorems are equivalent to \S 6 of [G] at least for
finite dimensional spaces.

On the other hand, we are
not aware of any other publications extending these to infinite
dimensional
spaces. Note that some of the best examples of measures of
symmetry --- e.g.\ the ``affine-invariant volume ratio", cf.\ [HC
p.~362]
--- use finite dimensional notions like volume, and thus
they cannot be extended to infinite dimensional spaces. Therefore
it is of some interest that the same Minkowski functional and
measure of symmetry can be used in that generality.

\proclaim {Theorem 6.5}.
Let $K\in \Kc$. Then $\al_K=0$ if and only if
$K$ is centrally symmetric. In this case $K^0=\{\x_0\}$ and the
symmetry center of $K$ is $\x_0$.

\proof
If $K$ is centrally symmetric with respect to a certain
point $\x_0$, then it is easy to see that $K^0=\{\x_0\}$ and of
course $\al_K=0$.

Assume now that $\al_K=0$. Then according to Lemma 6.3 we have
$K^0\ne\ey$ and $K^0=\{\x_0\}$ has exactly one element. But this
means that $\x_0\in L^0(K,\vb^*)$ for all $\vb^*\in S^*$, i.e.
$\x_0$ is a symmetry center for all $L(K,\vb^*)$ $(\vb^*\in S^*)$.
Since $K=K^1=\bcapi_{\vb^*\in S^*} L(K,\vb^*)$, and all
$L(K,\vb^*)$ are symmetric w.r.t. $\x_0$, we obtain that also $K$
is symmetric w.r.t. $\x_0$, as asserted.

\proclaim {Theorem 6.6}.
The mapping $f\colon \Kc( X )\to(0,1]$,
$f(K):= 1-\al_K$ is a {\it measure of symmetry\/}. That is, we
have
\begin{itemize}
\item[(i)]
$f(K)=1$ if and only if $K$ is centrally symmetric;
\item[(ii)]
$f$ is continuous (w.r.t. $\de$ on $\Kc)$;
\item[(iii)]
$f$ satisfies the ``superminimality condition"
$$
f(K_1+K_2)\ge \min \big\{f(K_1), \, f(K_2)\big\}\quad (K_1,K_2\in
\Kc).
$$
\end{itemize}

\proof (i) is equivalent to Theorem 6.5, and (ii) is a consequence
of Theorem 5.5. Finally, (iii) is the reformulation of our
Proposition 6.2 above.

\proclaim {Theorem 6.7}.
Let $K=M\ts N\sbt X=Y\ts Z$, $M\sbt
Y$, $N\sbt Z$, where $X,Y,Z$ are normed spaces and $K,M,N$ are
convex bodies. Then we have
\begin{itemize}
\item[(i)]
$K^\la=\{\x=(\yb,\zb)\colon \yb\in M^\la,\, \zb \in
N^\la\}=M^\la\ts N^\la$,
\item[(ii)]
$\al(K,\x)=\max\big\{\al(M,\yb),\,\al(N,\zb)\big\}
\quad \big(\x=(\yb,\zb)\big)$.
\end{itemize}

\proof
The two assertions are clearly equivalent, hence it
suffices to show (ii). Although a direct calculation is
possible, we take advantage of the proven equivalent forms of
$\al(K,\x)$.

Suppose first that $\y\in M$ and $\zb\in N$. Since homothetic
reduction can be written coordinatewise in $M$ and $N$,
definition (4.3) immediately provides
$$
\be(K,\x)=\be\big(M\ts N, (\yb,\zb)\big)=\min\big\{\be(M,\yb),
\be(N,\zb)\big\},
$$
which implies (ii) in view of Corollary 4.11 ii) and the fact that
${1-t\over 1+t}$ is decreasing.

Suppose now that $\x\notin K$, and apply (4.4). Obviously $K\cap
\big(\x+\la(K-\x)\big)=\ey$ if and only if either $M\cap\big(
\yb+\la(M-\yb)\big)=\ey$, or $N\cap \big(\zb+\la(N-\zb)\big)=\ey$.
First, for $\x\notin M$ and $\zb\notin N$ we find
$$
\vr(K,\x)=\vr\big(M\ts N,\,
(\yb,\zb)\big)=\max\big\{\vr(M,\yb),\,\vr(N,\zb)\big\}.
$$
Here now we refer to Corollary 4.11 iii) and the fact that the
function $\phi(t)={1+t\over 1-t}$ is increasing  $\big(t\in
(0,1)\big)$ to conclude (ii). Finally, if for example $\y\notin M$
and $\zb\in N$ we similarly have $\vr(K,\x)= \vr(M,\y)$ which by
$\al(N,\zb)\le 1 < \al(M,\y)$ yields (ii). The proof is complete.

\proclaim {Corollary 6.8}.
For $K=M\ts N\sbt X=Y\ts Z$, $M\sbt
Y$, $N\sbt Z$, where $X,Y,Z$ are normed spaces and $K,M,N$ are
convex bodies, we have
$$
\al_K=\max\{\al_M, \al_N\}.
$$

The special case of Corollary 6.8, where $\dim Z < \infty$, $M$ is
centrally symmetric but $N$ is not and $|N^{\al_N}|=1$, has been
mentioned in [HS].

\section{More on the sets $K^\la$}

While presenting the above results in Anogia-Crete in the summer
of 2001, Professor Rolf Schneider [Sch2] communicated us the following
nice formula which enables us to directly obtain some of the above
results as well.
(Note that this is also connected to a general theorem of
Schneider, c.f.\ [HC] p.~305 Theorem~2.)

\proclaim {Proposition 7.1 (R. Schneider)}.
For $\la\ge 1$ we have
$$
K^\la=\ol{{\la+1\over 2}K+{\la-1\over
2}(-K)}=\ol{(\la-1)C+K}.
\leqno(7.1)
$$

\proof
Let us check the last equality first. Note that by
convexity of $K$ we have $\mu K+\vr K=(\mu+\vr)K$ $(\forall
\mu,\vr\ge 0)$; hence
$$
{\la+1\over 2}K+{\la-1\over 2}(-K)=K+{\la-1\over
2}K+{\la-1\over 2}(-K)=K+(\la-1)C,
$$
as needed.

Next observe that the asserted formula (7.1) holds for all
centrally symmetric bodies (even if their center of symmetry is
not at $\ob$). In particular, for layers this was  already observed
in (3.1). Therefore, by using (3.1) and the fact that $L(\mu A+\nu
B,\vb^*)=\mu L(A,\vb^*)+\nu L(B,\vb^*)$ we get
$$
\begin{array}{rl}
L^\la(K,\vb^*)&= \dst{{\la+1\over 2}L(K,\vb^*)+{\la-1\over
2}\big(L(-K,\vb^*)\big)=}\\[4mm]
=& \dst{L\left({\la+1\over 2}K+{\la-1\over 2}(-K),
\vb^*\right)}.
\end{array}
\leqno(7.2)
$$
Thus applying the representation (2.9) now to
$\ol{{\la+1\over 2}K+{\la-1\over 2}(-K)} $ as well, we are led to
$$
\dsl{(7.3) \hfill K^\la= \bcapi_{\vb^*\in S^*}
L^\la(K,\vb^*)=\hfill\qquad
\cr\nomd
=\bcapi_{\vb^*\in
S^*} L\left({\la+1\over 2}K+{\la-1\over 2}(-K),
\vb^*\right)=\ol{{\la+1\over 2}K+{\la-1\over 2}(-K)}.
\cr}
$$

Note that for $X=\R^d$, compactness ensures closedness and there
is no need to take closure in formula (7.1).

\proclaim {Corollary 7.2}. For all $\la\ge 1$ and $\vb^* \in S^*$
we have the relations
\begin{itemize}
\item[i)]
$h(K^\la,\vb^*)= h\big( L^\la(K,\vb^*),\vb^*\big)= {\la+1\over
2} h(K,\vb^*)+ {\la-1\over 2} h(K,-\vb^*)$,
\item[ii)]
$w(K^\la,\vb^*)= w \big(L^\la (K,\vb^*),\vb^*\big)=
\la w (K,\vb^*)$,
\item[iii)]
$ L( K^\la,\vb^*) = L^\la (K,\vb^*)$.
\end{itemize}

Observe that the relation ``$\sbt$ " in (iii), or ``$\le$ " in the
first equality of (ii), are direct consequences of the definition
(2.12) even for arbitrary $\la\ge 0$.

\medskip

An analogous formula to (7.1) was pointed out to us, even for
$0\le \la \le 1$, by E. Makai [M]. To formulate this, one uses the
Minkowski difference
$$
A\sim B:=\{ \x\in X : \x+B\subset A\}.
\leqno(7.4)
$$

\proclaim{Proposition 7.3 (E. Makai)}.

For $0\le \la \le 1$ and $K \in \Kca$ we have
$$
K^{\la}=K \sim (1-\la)\ol{C}=\frac{1+\la}{2}K \sim
\frac{1-\la}{2}(-K).
\leqno(7.5)
$$

The key in proving Proposition 7.3  is a well-known lemma, see e.g.
formula (3.1.16) of [Sch], extended to infinite dimensional
spaces.

\proclaim{Lemma 7.4}.
For $K, M \in \Kca$ we have
$$
K\sim M=\bcapi_{\ub^*\in S^*}\{\x \in X : \,
<\ub^*,\x> \le h(K,\ub^*)-h(M,\ub^*)\}.
\leqno(7.6)
$$

\proof
Just follow the argument in [Sch p. 134], where the restriction to
finite
dimensional spaces is of no relevance.

{\bf Proof of Proposition 7.3.} By (7.6) and (2.4) we have
$$
\begin{array}{rl}
K\sim & (1-\la)\ol{C} =
\\
=& \bcapi_{\vb^*\in S^*}\{\x \in X : \quad <\vb^*,\x> \le
h(K,\vb^*)-\frac{1-\la}{2}w(K,\vb^*)\} =
\\
=& \bcapi_{\vb^*\in S^*}\{\x \in X : \quad
<\vb^*,\x> \le
\frac{1+\la}{2}h(K,\vb^*)-\frac{1-\la}{2}h(K,-\vb^*)\}.
\end{array}
\leqno(7.7)
$$

Grouping together opposite pairs $\vb^*, -\vb^* \in S^*$ and
reflecting back to (2.10)--(2.12) yields the first equation of
(7.5). Since $h(-K,\vb^*)=h(K,-\vb^*)$, it is even easier
to see that
$$
\dsl{
(7.8) \qquad \qquad
\frac{1+\la}{2}K \sim \frac{1-\la}{2}(-K)=
\hfill\quad
\cr\nomd
= \bcapi_{\vb^*\in S^*} \left\{\x \in X : \quad
<\vb^*,\x> \le
\frac{1+\la}{2}h(K,\vb^*)-\frac{1-\la}{2}h(K,-\vb^*)\right\}\,.
\cr\nomd}
$$
This proves the second equation of (7.5).

\medskip

At this stage we point out that the ``level sets'' $K^\la$ were known
in geometry under an equivalent definition. In fact, for any
convex body $K\in \Kca$ Hammer [H] has introduced the associated
bodies

$$
C(K,\rho):=\left\{\begin{array}{ll}
\bcapi_{\x\in \pa K}\{\x + \rho (K-\x)\} \qquad , & 0\le\rho\le 1 \\
\overline{\bcupi_{\x\in \pa K}\{\x + \rho (K-\x)\}} \qquad, &
\rho\ge 1 \quad
\end{array}\right. ,
\leqno(7.9)
$$
naturally without the closure, since Hammer has considered only
finite dimension.

Hammer's definition (7.9) is closely related to the functional
$\si(K,\x)$ of (4.7); in fact, Hammer has defined these associated
bodies as level
sets of the functional $\si(K,\x)$ in Minkowski's measure of
symmetry,
c.f. [G, p. 246-247].
Thus, one can expect that there is an equivalence, like for the
functionals, also even for the level sets as well. Using Corollary
4.11 (the
relations between $\al(K,\x)$ and $\si(K,\x)$),
it is not hard to obtain this equivalence for $\dim X < \infty$.
In fact, Hammer's Theorem 5 can be regarded as a predecessor to
our
Corollary 4.11 with several versions of the equivalent extremal
ratios
mentioned there.
However, there is another way to get the equivalence.
In [Sch p. 141] one finds the finite dimensional version of

$$
C(K,\rho)=\left\{\begin{array}{ll}
K \sim 2(1-\rho)\CK  \qquad , & 0\le\rho\le 1 \\
\overline{K+2(\rho-1)\CK}  \qquad , & \rho\ge 1\quad
\end{array}\right.
\leqno(7.10)
$$
without the closure; the extension to $\dim(X)=\infty$ is
standard.
But comparing (7.10) with Propositions 7.1 and 7.3, we are led to
$$
C(K,\rho)=K^{2\rho-1} .
\leqno(7.11)
$$
\medskip

Let us continue with a general statement, obtained in a different
way by H\"ormander [H\"o2] in 1955.

\proclaim {Lemma 7.5 (H\"ormander)}.
If $K$ and $K'$ are arbitrary non-empty bounded
convex sets in $X$, their Hausdorff distance satisfies
$$
\delta(K,K')=\sup_{\vb^*\in S^*}|h(K,\vb^*)-h(K',\vb^*)|.
\leqno(7.12)
$$

\proof Recall that already in the proof of Lemma 5.4 we have
mentioned the `` $\ge$ " direction of the stated formula as a trivial
consequence of definition (2.14) of the Hausdorff distance. To
prove the other direction, it suffices to restrict ourselves to
closed convex sets $K$ and $K'$ since by taking the closure does not
change either sides of the stated inequality.

First, let's assume that $\x\in K$ and
$B(\x,r)\cap K' =\ey$, with $r>0$. Since both $B(\x,r)$ and $K'$
are
convex, and $\intt B(\x,r)\ne\ey$, by the separation theorem of
convex sets
(c.f.\ [Ru], Theorem 3.4) we can find a $\vb^*\in S^*$ such that
$$
h(K',\vb^*)=\sup_{K'}\vb^*\le \inf_{B(\x,r)}\vb^*=\<\vb^*,\x\>-r
\le \sup_K \vb^* -r=h(K,\vb^*)-r.
\leqno(7.13)
$$
Second, let us consider the case when $\y\in K'$ and $B(\y,R)\cap
K = \ey$,
where $R>0$.
Working as above, for some appropriate $\ub^*\in S^*$ we have
$$
h(K,\ub^*) \le h(K',\ub^*)-R.
\leqno(7.14)
$$
By (2.14) we have for any $\epsilon > 0$ either some $\x\in K$ and
$r>\delta - \epsilon$, or some $\y \in K'$ and $R>\delta -
\epsilon$ with the
above properties; hence either $h(K',\vb^*)\le h(K,\vb^*)+\epsilon
- \delta$,
or $h(K,\ub^*)\le h(K',\ub^*)+\epsilon-\delta$. In either case we
conclude
$\sup_{S^*}|h(K,\cdot)-h(K',\cdot)|\ge \delta$, as needed.

As a first application, we present a sharp form of Corollary 5.9.

\proclaim {Corollary 7.6}.
For $\la \ge 1$ we have
$$
\de\big(K^\la,\la C\big)\le D(K)-{1\over 2}w(K).
\leqno(7.15)
$$

\proof Since $\la \ge 1$, we can apply Schneider's Formula, or its
Corollary 7.2 (i) which is more suitable to us here. For an
arbitrary $\vb^*\in S^*$ and by using (2.4) we obtain
$$
\dsl{
h(K^\la,\vb^*)-h(\la C,\vb^*) =
{\la+1\over 2} h(K,\vb^*)+ {\la-1\over 2} h(K,-\vb^*) -\la
h(C,\vb^*)
\hfill\quad
\cr\nomd
= {1\over 2}(h(K,\vb^*)- h(K,-\vb^*)) \le h(K,\vb^*) - {1 \over
2}w(K)
\le D(K) - {1 \over 2}w(K).
\cr\nomd}
$$
An analogous calculation gives $h(\la C,\vb^*)-h(K^\la,\vb^*)\le
D(K)-\frac12 w(K)$. Since by Lemma 7.5 we have $\de\big(K^\la,\la
C\big)=\sup_{S^*} |h(K^\la,\cdot)-h(\la C,\cdot)|$, the proof of
the Corollary follows.

\proclaim {Remarks 7.7}. \rm (i) Corollary 7.6 is quite sharp.
For instance, let $K$ be a ball $B(\ab,r)$ of radius $r$. Then $w(K)=
d(K)=2r$, $C=C=B(\ob,r)$, $K=C+\ab$ and $D(K)=a+r$, where
$a=\|\ab\|$. For this $K$ and for all $\la\ge 0$ we have
$$
\de\big(K^\la, \la C\big) =a = D(K)-{1\over 2}w(K).
$$

(ii) Remarkably, (7.15) does not hold for $\la < 1$. To construct
an example, consider in $X=\R^2$ the triangle $K:={\rm
conv}\{(10,10);(16,10);(10,16)\}$. In this case $w(K)=3\sqrt{2}$ and
$D(K)=\sqrt{10^2+16^2}=18.86...$, while
$D(K)-{1\over2}w(K)=16.74...$. On the other hand
$$
C(K)={\rm conv}\{(3,0);(0,3);(-3,3);(-3,0);(0,-3);(3,-3)\},
$$
and for $\la={1\over3}$ we have $K^\la=\{(12,12)\}$. By taking
$\y:=(0,-1)\in \la C$, we can calculate that ${\rm
dist}(\y,K^\la)=\sqrt{12^2+13^2}=17.69...$ which exceeds the above
value of $D(K)-{1\over2}w(K)=16.74...$. Although all computed
examples suggest that $D(K)$ should be a generally valid upper
bound for all $\la > \al_K$, we are not aware of a sharp bound on
$\delta (K^\la,\la C)$ in the case $\la <1$.

\proclaim {Corollary 7.8}.
Let $1\le \la,\mu$ be arbitrary. Then we have
$$
{(K^\la)}^\mu=K^{\la\mu}.
\leqno(7.16)
$$

\proof
$$
\dsl{{(K^\la)}^\mu=\ol{(\mu-1)C(K^\la)+K^\la}=
\cr\nomd
=\ol{(\mu-1)C\big(\ol{(\la-1)C+K}\big)+
\ol{(\la-1)C (K)+K}}=
\cr\nomd
=\ol{(\mu-1)\big\{(\la-1)C
(K)+C\big\}+(\la-1)C+K}=
\cr\nomd
=\ol{(\mu\la-1)C+K}=K^{\la\mu}.
\cr}
$$

\proclaim {Remark 7.9}.
\rm Note that ${(K^\la)}^\mu \sbt
K^{\la\mu}$ $(\la>0,\,\mu\ge 1)$, but not for all $\la,\mu>0$.

\medskip

In fact, a further generalization of the above can be found in
[Sch] where the relative parallel bodies of $K$ with respect to
$D\in \Kca$ are defined as
$$
K_{\rho,D}:=\left\{ \begin{array}{ll}
\ol{K+\rho D} \qquad \quad\,\,\,\, ,& \rho \ge 0
\\
K \sim (-\rho D) \qquad ,& \rho <0 \quad.
\end{array}\right.
\leqno(7.17)
$$
Here we use closure in (7.17) only in view of the case $\dim X =
\infty$. Clearly, $K_{\rho,D}$ is a convex body for $\rho >
\rho_0$, with some $\rho_0>0$ and, in analogy to \S 2, basic
properties of this parametric set sequence can be easily seen.

\proclaim{Lemma 7.10}.
For $K,M,D \in \Kca$ and $\rho, \sigma \in \R$ we have
$$
K_{\rho,D}+M_{\si,D} \subset \ol{(K+M)_{\rho+\si,D}}.
\leqno(7.18)
$$

\proof Taking into consideration (7.6), the obvious analogy for the closure of
the Minkowski sum, and by using the fact that the  convexity of
$K,M$ and $D$ implies that $K_{\rho,D} + M_{\si,D}$ is also convex,
we obtain
$$
\begin{array}{rl}
K_{\rho,D}+M_{\si,D} \subset
\qquad \qquad \qquad \qquad \qquad \qquad \qquad \qquad \qquad
\qquad
\\
\subset \bcapi_{\vb^* \in S^*}
\{ X_{h(K,\vb^*)+\rho h(D,\vb^*)}(\vb^*) + X_{h(M,\vb^*)+\si
h(D,\vb^*)}(\vb^*)\}=
\\
= \bcapi_{\vb^* \in S^*}
X_{h(K+M,\vb^*)+(\rho+\si)h(D,\vb^*)}(\vb^*) =
(\ol{K+M})_{\rho+\si,D}.\quad
\end{array}
\leqno(7.19)
$$
This approach gives rise to the more general definition
$$
\al _D (K,\x):= \inf \{ \rho \in \R : \x \in K_{\rho, D} \}.
\leqno(7.20)
$$
Note that with this extension $\al (K,\x) = 1 + \al
_{\ol{C}}(K,\x)$. Therefore, a number of properties of $K^\la$ and
$\al (K,\x)$ can be obtained from more general statements for
$K_{\rho,D}$ and $\al_D(K,\x)$. However, we don't need them for
our present applications and we leave this subject now.

\section{Small values of $\al$ and the critical set of $K$}

Here is an interesting question in geometry: what kind of
quasi-centre can be found for a nonsymmetric convex $K$? A usual
choice is the centroid, but there are other important
definitions like e.g.\ the {\it Santal\'o point\/} (c.f.\ [HC],
p.~165) or the {\it analytic center\/} (c.f.\ [HC], p.~661) of a
convex body.

The general notion of a center of a convex body with respect to
the {\it center function\/} $\tr\colon \Kc^d\ts \R^d\to \R_+$ was
introduced by Huard in 1967, see [Hu] or [HC p.~660]. According to
Huard, a general center function has to satisfy

$$
\begin{array}{ll}
\hbox{ (i) } & \tr(K,\x')=\min\limits_{\x\in
K}\tr(K,\x)\quad (\forall \x'\in\pa K,\,K\in \Kc^d)\\[4mm]
\hbox{ (ii) } & \tr(K,\x')>\min\limits_{\x\in K} \tr(K,\x)\quad
(\forall \x'\in\intt K,\,K\in \Kc^d)\\[4mm]
\hbox{ (iii) } & \tr(K,\x)\ge \tr(K',\x) \quad (\forall \x\in
K'\sbt K, \, K,K'\in \Kc^d)\\[4mm]
\hbox{ (iv) } & \max\limits_{\x\in K}\tr (K,\x)=\tr (K,\x^*)
\hbox{
for a unique point }\x^*\in K \,(K\in \Kc^d).
\end{array}
\leqno(8.1)
$$

Webster (see [W], p.~319) presents another approach. Starting from
the observation that the symmetry center of a centrally symmetric
convex body bisects every chord of the body through the centre, Webster
is looking for a point in $K$ which divides all chords
through it ``relatively evenly".
Theorem 7.1.5 in [W] states that for all $K\in
\Kc^d$ there exists a point $\cb=\cb_K\in\intt K$, such that for
any chord $l\cap K=[\ab,\bb]\ni \cb$, the segments always satisfy
${\|\bb-\cb\|\over \|\ab-\bb\|}\le {d\over d+1}$.
In other words, we have $\om(K,\cb)\le {d\over d+1}$.
This translates to the inequality $\al(K,\cb)\le {d-1\over d+1}$,
see Corollary 4.11 (ii).

Hence, Webster's Theorem 7.1.5 is seen to be equivalent to the
assertion that $\al_K\le {d-1\over d+1}$ and that $K^{{d-1\over
d+1}}\ne\ey$. This leads to the natural idea of considering
$K^{\al_K}$
(or, equivalently the function $\tr(K,\x):= 1- \al(K,\x)$) as a
candidate
for a center function. In fact, by an equivalent definition, the
same set was
already known in geometry for $\dim X < \infty$ as the
{\it critical set} of $K$, see [H] or [Sch, p. 141].

In [KR] the question of cardinality of $K^{\al_K}$ was raised.
However, it
turns out that the answer preceded this question by half a
century.
It is true that $|K^{\al_K}|=1$ for dimension $d=2$, see Neumann
[N], but in
general Sobczyk, see [H, p.792] and [HS], has answered this question
in the negative.

\proclaim {Example 8.1 (Sobczyk)}.
\rm Let $K=\tr\ts I=\big\{(x,y,z)\in
\R^3\colon 0\le x,y\le x+y\le 1,$
$-1\le z\le 1\big\}\sbt \R^3$. Then
$\al_K={1\over 3}$ and $K^{1/3}=\left\{\left({2\over 3}, {2\over
3},z\right)\colon |z|\le {1\over 3}\right\}$.

Indeed, for the triangle $\tr=\big\{(x,y)\colon 0\le x,y\le x+y\le
1\big\}$ we have $\al_{\tr}={1\over 3}$, while $\tr^{1/3}=
\big\{({2\over3},
{2\over3})\big\}$. On the other hand $\al_I=0$ and $I^{1/3}
=\left[-{1\over
3},{1\over 3}\right]$. Hence, by Theorem 6.7 (i)
$K^{1/3}=\left\{\left({2\over 3},{2\over 3}\right)\right\}\ts
\left[-{1\over 3},{1\over 3}\right]$ and by Corollary 6.8
$\al_K=1/3$.

\proclaim{Remarks 8.2}
\rm
(i) Observe that in the general context the above construction
is not a
``center function" in the sense of Huard, as property (iv) is seen
to fail by Example 8.1.
\newline \indent \rm (ii)
However, the most important fact is that (iii) of (8.1)
does not hold. To see this just consider a large cube $K$ and a small one
$K'$ in one of its corners. While for $K'$ the generalized Minkowski
functional  $\al(K',\x)$ assumes all $\le 1$ values on $K'$,
it is clear that only values close to 1 will be assumed by
$\al(K,\x)$ for $\x\in K'$.
\newline \indent \rm (iii)
Still, the size and properties of the critical set $K^{\al_K}$
are of interest. For the finite dimensional case Klee [K] has
obtained a rather
precise description of the critical set. For finite dimensional
spaces  Klee's estimate has as follows
$$
{r_K\over 1-r_K}+\dim K^{\al_K} \le \dim X,
$$
where his $r_K$ arises by using $\omega(K,\x)$ and so by
Corollary 4.11 one
obtains $r_K:=\inf \omega(K,\dot)={1+\al_K \over 2}$.
This is equivalent to
$$
\frac{1+\al_K}{1-\al_K} \le \codim K^{\al_K}.
$$
In this form the assertion can be interpreted even if $\dim
X=\infty$.
However, the whole argument of [K] relies essentially on
compactness, so
much so that even the objects he starts with do not exist in
general.
Hence, not only a description of the critical set is missing but even
Question 6.4 seems to be complicated for $\dim X = \infty$.

\proclaim {Theorem 8.3 (Minkowski, Radon)}.
For any $K\sbt \R^d$, $\al_K\le
{d-1\over d+1}$. Moreover, we always have $\w_K\in
K^{{d-1\over d+1}}$ for the centroid $\w_K$ of $K$, and
$\al(K,\w_K)={d-1\over d+1}$ only if $K$ is a finite cone.

This result is well-known and was obtained already by Minkowski
[Mi] for $d=2,3$
and Radon [R] for $d\in \N$. See also \S 6.2 (i) of [G] and the
references
therein. The reader can easily observe that Theorem 8.3
contains also the above mentioned Theorem 7.1.5 of [W].

\proclaim {Example 8.4}.
\rm Consider any simplex $S\sbt \R^d$. It is
easy to see that $S^{{d-1\over d+1}}=\{\w_S\}$ and thus
$\al(S,\w_S)=\al_S={d-1\over d+1}$.

It is also known that in $\R^d$ simplices are the only extremal
bodies with
$\al_K={d-1 \over d+1}$, see (4.6) of [K] or \S 6.2 (i) in [G].

\proclaim {Example 8.5 (Hammer)}.
\rm Let $d=2$, and $K=\big\{(x,y)\in
\R^2\colon 0\le y$, $x^2+y^2\le 1\big\}$ be the upper half of
the unit disc. Then $\w_K=\left(0, {4\over 3\pi}\right)$, $\al_K=3
-2\sqrt 2<{1\over 3}$, and $K^{\al_K}
=\big\{(0, \sqrt 2-1)\big\}$. Thus $\w_K\notin K^{\al_K}$.

With respect to Theorem 8.3, it is worth noting that for vector
spaces of infinite dimension no general upper bound, lower than
1, can be given to $\al_K$ of $K\in \Kc( X )$.
Compare also to Corollary 3.8 (vi).

\proclaim {Proposition 8.6}.
Let $X$ be an infinite dimensional space. Then
$$
\sup\big\{\al_K\colon K\in \Kc(X)\big\}=1.
$$

\proof Let $d\in \N$ be arbitrary and let $Y$ be any
$d$-dimensional subspace of $X$. Suppose $\x_1,\dots,\x_d$ are
linearly independent vectors of $Y$ and let $Z$ be a complement to
$Y$ in $X$. That is $X=Y\oplus Z$, where $Z$ is a closed
subspace such that $Y\cap Z=\{\ob\}$, $Y+Z=\{\yb+\zb\colon \yb\in
Y,\,\zb\in Z\}=X$. The existence of a complement to each finite
dimensional subspace of a normed space $ X $
is a well known result in functional analysis, see for instance
[Ru, Lemma 4.21]

Now consider the simplex
$S=\con\{\ob,\x_1, \dots, \x_d\}$ in $Y$ and
define $K:= S\oplus B_Z$, where $B_Z$ is the unit ball  of $Z$.
This is the direct sum construction in Corollary 6.8. So,
by Example 8.4 it can be easily seen that
$\al_K=\max\{\al_S,\al_{B_Z}\}=\al_S={d-1\over d+1}
\le \sup_{\Kc( X )}\al_K$ holds for all $d\in \N$,
proving the assertion.

\section{Applications in approximation theory}

To illustrate the use of the above generalized Minkowski
functional, we describe now two applications from the field of
multivariate polynomial extremal problems. Our general setting
here again will require that $K\sbt X$ is a convex body, but we
also introduce the space of continuous polynomials of degree at
most $n$, denoted by $\Pc_n(X)=\Pc_n$, from $X$ to $\R$. In
infinite dimensional normed spaces, the space of continuous (or,
equivalently, bounded) polynomials depends on the norm too. That
is, for different, non-equivalent norms, even the space of
polynomials (and in particular the space of continuous linear
functionals) will be different. However, once $K$ is a fixed
convex {\it body}, that is, bounded and also containing a
ball, we conclude that all such norms are equivalent (and actually
equivalent to the norm generated by the central symmetrization
$C$ of $K$ by means of its Minkowski functional $\vp_{C}$).
Since we always normalize the polynomials defined on $K$ with
respect to the sup-norm, that is
$$
{\|p\|}_K:= \sup_K|p|\qquad \big(p\in \Pc_n(X)\big),
\leqno (9.1)
$$
it is clear that the resulting set of polynomials $\big\{p\in
\Pc_n(X),\;{\|p\|} _K\le 1\big\}$ will always be the same,
provided $K$
is fixed.

To define, in general, polynomials we can follow e.g.\ H\"ormander
[H\"o],
Appendix A. We say that {\em $p$ is a polynomial of degree at most
$n$} if for arbitrary fixed $x,y\in X$ the function $p(x+ty)\colon
\R\to \R$
is a real algebraic polynomial in $t\in \R$ of degree at most $n$.

A polynomial $p\in \Pc_n(X)$ is homogeneous of degree $m$, if
it satisfies $p(t\x)=t^m p(\x)$ $(\forall t\in \R,\,\forall \x\in
X)$. It is easy to see that every polynomial $p\in \Pc_n(X)$ can
be written as
$$
p=\sum_{k=0}^n p_k\qquad (p_k \hbox{ is homogeneous of degree
}k).
\leqno(9.2)
$$
Note that in the previous expression the ``leading term" $p_n$ of $p$
is a homogeneous polynomial of degree $n$.

Continuous polynomials are Fr\'echet differentiable everywhere
and with any order of differentiation. In particular,
the gradient of $p$ at a point $\x\in X$ is
$$
\grad p(\x)\colon X\to\R, \qquad \grad p(\x)\in X^*.
\leqno(9.3)
$$
A pointwise estimate of $\grad p(\x)$, for $\x\in K$, under
normalization of the norm defined by (9.1) is the so-called
{\em Bernstein problem} or {\em Bernstein type inequality}.
For centrally symmetric convex bodies we have the following
generalization of the classical one century old
result of S. Bernstein.

\proclaim {Theorem (Sarantopoulos, 1991 [S])}. Let $K$ be $\ob
-$symmetric , $\x\in\intt (K)$ and $p\in \Pc_n$. Then we have
$$
\big|\grad p(\x)\big|\le {2n{\|p\|}_K\over
w(K)\sqrt{1-\vp_K^2(\x)}}.
\leqno(9.4)
$$

Here $\vp_K$ is the Minkowski functional of the $\ob-$symmetric
body $K$. Note that actually $(9.4)$ was proved by Sarantopoulos
in the equivalent setting when $K$ is the unit ball, i.e.\ when
$\|\cdot\|=\vp_K$. Observe that in this case $w(K)=2$. However,
the two formulations are equivalent, so for further use we choose
the above version.

For more general, not necessarily symmetric sets, various
estimates have been proved, see [A], [B], [S1], [S2] and [S3]. All
involve some geometric techniques to quantify the location of $\x$
within $K$, i.e.\ the respective distance of $\x$ from the
boundary. These ``Bernstein factors" are sometimes quite
complicated, and even the comparison of the various results is
rather difficult. We quote here a result for the finite
dimensional case.

\proclaim {Theorem (Kro\'o--R\'ev\'esz, 1999 [KR])}.
Let $K\sbt
\R^d$ be a convex body, $\x\in\intt K$, and $p\in \Pc_n$ be
arbitrary. Then we have
$$
\big|\grad p(\x)\big|\le {2n{\|p\|}_K\over
w(K)\sqrt{1-\al(K,\x)}}.
\leqno(9.5)
$$

In [KR] the possibility of extending the above result to arbitrary
Banach spaces was mentioned, but the details were not
worked out.

Note the similarity of (9.5) to (9.4), but $\al$ is
not raised to the power 2. However, it is our conjecture that $\al$
should be raised to the power 2.

\proclaim Conjecture. If $K\sbt X$ is a convex body and if $\x\in\intt
K$ and $p\in \Pc_n$ are arbitrary, then we have
$$
\big|\grad p(\x)\big|\le {2n{\|p\|}_K\over
w(K)\sqrt{1-\al{(K,x)}^2}}.
\leqno(9.6)
$$

Let us turn our attention to another group of problems, the so
called Chebyshev-type extremal problems concerning growth of
real polynomials. The question has as follows: `` How large can be a
polynomial at a point $\x\in X$, or when $\x\to\iy$? " More
precisely, we are interested in determining for arbitrary
$\vb\in X $ (say, with $\|\vb\|=1$)
$$
A_n(K,\vb):= \sup\big\{p_n(\vb)\colon p\in \Pc_n\hbox{ satisfying
(9.2), } {\|p\|}_K\le 1\big\},
\leqno(9.7)
$$
or for some fixed $\x\in X$
$$
C_n(K,\x):= \sup\big\{p(\x):p\in \Pc_n,\ {\|p\|}_K\le 1\big\}.
\leqno(9.8)
$$
Note the appearance of the $n$-homogeneous part $p_n$ in (9.7).
Both problems are classical and fundamental in the theory of
approximation, see e.g.\ [MMR] or [Ri] for the one and a half
century old single variable result and its many consequences,
variations and extensions.

As we have mentioned in the introduction, even the above formulation
(2.12)-(2.13)
of the definition of $\al(K,\x)$ was applied first in work on
these
questions, particularly on (9.8), where a quantification of the
position of
$\x$ with respect to $K$ is needed. After the results in $\R^d$ by
Rivlin--Shapiro [RiSh] and Kro\'o--Schmidt [KS], we were able to settle
the general case in a satisfactory way. We have the following result.

\proclaim {Theorem (R\'ev\'esz--Sarantopoulos, 2001, [RS])}.
If $K\sbt X$ is an arbitrary convex body and if $\x\in X\sm K$
is arbitrary, then we have
$$
C_n(K,\x)=T_n\big(\al(K,\x)\big).
\leqno(9.9)
$$
Moreover, $C_n(K,\x)$ is actually a maximum, attained by
$$
P(\x):= T_n\big(t(K,\vb^*,\x)\big).
\leqno(9.10)
$$
Here $T_n$ is the classical Chebyshev polynomial
$$
\begin{array}{rl} T_n(x)&:= 2^{n-1}\prod\limits_{j=1}^n
\left(x-\cos\left({(2j-1)\pi\over
2n}\right)\right)=\\[4mm]
&= {1\over 2}\left\{{(x+\sqrt{x^2-1})}^n
+{(x-\sqrt{x^2-1})}^n\right\},
\end{array}
\leqno(9.11)
$$
$t(K,\vb^*,\x)$ is the linear expression defined in $(4.1)$, and
$\vb^*$ is some appropriately chosen linear functional from $S^*$.

Note that actually the restriction $\x\notin K$ is natural, as
$p\ev 1\in \Pc_n$, and thus for $\x\in K$ we always have
$C_n (K,\x)=1$. Here we can observe that our results on the growth
of $\al(K,\x)$ together with (9.9) give strong indications even
for the other Chebyshev problem, as we know that the Chebyshev
polynomial itself has leading coefficient $2^{n-1}$. Indeed, we have
the following result.

\proclaim {Theorem (R\'ev\'esz--Sarantopoulos, 2001 [RS])}.
Let $K\sbt X$ be an arbitrary convex body and let $\vb\in X$. Then
we have
$$
A_n(K,\vb)={2^{2n-1}\over \tau{(K,\vb)}^n},
\leqno(9.12)
$$
and the supremum is actually a maximum attained by a polynomial
of the form $(9.10)$ with some appropriately chosen $\vb^*$
satisfying (2.25).

Based on the determination of these extremal quantities, other
related questions were already addressed in approximation theory,
such as the\break uniqueness of the extremal polynomials, or
the existence
of the so-called universal majorant polynomials. These, in turn,
have consequences e.g.\ concerning the approximation of convex
bodies by convex hulls of algebraic surfaces. For further
details we refer to [Kr] and [R\'e].

\section{Acknowledgement}

The authors would like to express their indebtedness to Imre
B\'ar\'any,
Apostolos Giannopoulos, Endre Makai and Rolf Schneider for many
useful
references, suggestions and comments.

\section{References}

{\parindent=40pt

\begin{itemize}

\item[{[A]}]
Andrianov, A. V., On some open problems for algebraic
polynomials on bounded convex bodies, {\it East J. Approx.,\/}
{\bf 5,
No 1} (1999), 117--123.

\item[{[B]}]
Baran, M., Bernstein type theorems for compact sets in $\R^n$
revisited, {\it J. Approx. Theory,\/} {\bf 79} (1994), 190--198.

\item[{[B\'a]}]
B\'ar\'any, I., The technique of $M$-regions and cap
coverings: a survey, {\it manuscript}.

\item[{[BL]}]
B\'ar\'any, I.--Larman, D.G., Convex bodies,
economic cap coverings, random polytopes, {\it Mathematika\/} {\bf
35}
(1988), 274--291.

\item[{[D]}]
Dineen, S., Complex analysis on infinite dimensional
spaces, {\it Springer Monographs in Mathematics,\/} Springer
Verlag, 1999.

\item[{[E]}]
Edwards, R. E., Functional Analysis,
Holt--Rinehart--Winston, New York--Toronto--London, 1965.

\item[{[G]}]
Gr\"unbaum, B., Measures of symmetry for convex sets, {\it
Proceedings
of Symposia in Pure Mathematics}, {\bf VII} AMS, 1963, 233--270.

\item[{[H]}]
Hammer, P. C., Covex bodies, associated with a convex body, {\it
Proc. Amer. Math. Soc.} {\bf 2} (1951), pp. 781--793.

\item[{[HS]}]
Hammer, P. C. -- Sobczyk, A., Critical points of a convex body,
Abstract 112,
{\it Bull. Amer. Math. Soc.} {\bf 2} (1951), 127.

\item[{[HC]}]
Handbook of Convex Geometry, Vol A,
{\it P. M. Gruber $\&$ J. M. Wills, eds.,\/},
North Holland, Amsterdam--London--New York--Tokyo, 1993.

\item[{[H\"o]}]
H\"ormander, L., Notions of Convexity,
{\it Progress in Mathematics Vol.\/} {\bf 127},
Birkhauser, Basel--Boston-- Berlin, 1994.

\item[{[H\"o2]}]
H\"ormander, L., Sur la fonction d'appui des ensembles convexes
dans un espace
localement convexe, Arkiv f\"or Mat. {\bf 3} (1955), 181-186.

\item[{[Hu]}]
Huard, P., Resolution of mathematical programming
with nonlinear constraints by the method of centers, {\it in
Nonlinear Programming (ed. J. Aladie),\/} North--Holland, 1967,
209--219.

\item[{[K]}]
Klee, V. L., The critical set of a convex body, {\it Amer. J.
Math.},
{\bf 75} (1953), 178--188.

\item[{[Kr]}]
Kro\'o, A., Universal polynomial majorants on convex bodies, {\it
J. Approx. Theory} {\bf 111} (2001), 220-232.

\item[{[KR]}]
Kro\'o, A.--R\'ev\'esz, Sz. Gy., On Bernstein and
Markov-type inequalities for multivariate polynomials on convex
bodies, {\it J. Approx. Theory,\/} {\bf 99} (1999), 134--152.

\item[{[KS]}]
Kro\'o, A.--Schmidt, D., Some extremal problems for
multivariate polynomials on convex bodies, {\it J. Approx.
Theory\/}
{\bf 90/3} (1997), 415--434.

\item[{[M]}]
Makai, E., Personal communication, July 2002.

\item[{[MMR]}]
Milovanovi\'c, G. V.--Mitrinovi\'c, D. S.--Rassias, Th.
M.,Topics in Polynomials: Extremal Problems, Inequalities,
Zero, World Scientific, Singapore, 1994.

\item[{[Mi]}]
Minkowski, Allgemeine Lehrs\"atze \"uber konvexe Polyeder,
{\it Nachr. Ges. Wiss. G\"ottingen}, 1897, 198-219. (=Ges. Abh.
vol. {\bf 2}
pp. 103-121, Leipzig--Berlin, 1911.)

\item[{[N]}]
Neumann, B. H., On some affine invariants of closed convex
regions,
{\it J. London Math. Soc.} {\bf 14} (1939), 262-272.

\item[{[P]}]
Phelps, R. R., Convex Functions, Monotone Operators
and Differentiability, {\it Lecture Notes in
Mathematics\/} \# {\bf 1364}, {\it 2nd
ed.,\/} Springer Verlag, Berlin--Heidelberg--New York--Tokyo,
1993.

\item[{[R]}] Radon, J., \"Uber eine Erweiterung des Begriffs der konvexen
Funktionen, mit einer Anwendung auf die Theorie der konvexen
K\"orper,
{\it S.-B. Akad Wiss. Wien} {\bf 125} (1916). 241-258.

\item[{[R\'e]}]
R\'ev\'esz, Sz., Gy., Uniqueness of multivariate Chebyshev-type
extremal\break polynomials for convex bodies, {\it East J.
Approx.},
{\bf 7} (2001), 205-240.

\item[{[RS]}]
R\'ev\'esz, Sz. Gy.--Sarantopoulos, Y., Chebyshev's
extremal problems of polynomial growth in real normed spaces,
{\it J. of Contemporary Analysis and Applications,\/},
{\bf 36} \#5 (2001), 62-81.

\item[{[Ri]}]
Rivlin, T. J., The Chebyshev Polynomials, Wiley, New
York--Sydney--London--Toronto, 1974.

\item[{[RiSh]}]
Rivlin, T. J.--Shapiro, H. S., A unified approach
to certain problems of approximation and minimization, {\it J.
Soc.
Ind. Appl. Math.\/} {\bf 9} (1961), 670--699.

\item[{[Ro]}]
Rockafellar, R. T., Convex Analysis, Princeton University Press,
\break Princeton, 1970.

\item[{[Ru]}]
Rudin, W., Functional Analysis (2nd ed.), McGraw- Hill, Inc., New
York, 1991.

\item[{[S]}]
Sarantopoulos, Y.,
Bounds on the derivatives of polynomials on Banach spaces,
{\it Math. Proc. Cambridge Philos. Soc.\/}
{\bf 110} (1991), 307--312.

\item[{[Sch]}]
Schneider, R., Convex Bodies: the Brunn--Minkowski theory,
Encyclopedia in
Mathematics and Applications, vol. {\bf 44}, Cambridge University
Press,
1993.

\item[{[Sch2]}]
Schneider, R., {\it Personal communication,\/} August 2001.

\item[{[S1]}]
Skalyga, V. I., Analogues of the Markov and
Bernstein inequalities on convex bodies in Banach spaces
{\it Izvestiya: Mathematics\/} {\bf 62:2}(1998), 375--397.

\item[{[S2]}]
Skalyga, V. I., Analogs of the Markov and Schaeffer--Duffin
Inequalities for Convex Bodies, {\it Mathematical Notes,\/} {\bf
68:1} (2000), 130--134.

\item[{[S3]}]
Skalyga, V. I., Estimates of the Derivatives of
Polynomials on Convex Bodies, {\it Proc. Steklov Inst. Math.\/}
{\bf
218} (1997), 372--383.

\item[{[W]}]
Webster, R., Convexity, Oxford University Press,
Oxford, 1994.

\item[{[We]}]
Werner, E., Illumination bodies and the affine
surface area, {\it Studia Math.\/} {\bf 110} (1994), 257--269.

\end{itemize}

}

\noindent\begin{minipage}{0.5\hsize}
{A. R\'enyi Institute of Mathematics\\
Hungarian Academy of Sciences,
\\Budapest, P.O.B. 127, 1364 Hungary\\ }
{email: revesz@renyi.hu}\\
\end{minipage}
\begin{minipage}{0.5\hsize}
{National Technical University\\
School of Applied Math.Phys. Sci.,\\
Department of Mathematics,\\
Zografou, 15780 Athens, Greece\\ } {email: ysarant@math.ntua.gr}
\end{minipage}
\end{document}